
\input amstex

\documentstyle{amsppt}

\magnification=\magstep1

\def\ds{\displaystyle}
\def\wt{\widetilde}

\def\dom{\text{dom}}

\NoBlackBoxes

\input BoxedEPS
\SetTexturesEPSFSpecial
\HideDisplacementBoxes


\topmatter
\title Calculus on the Sierpinski gasket I: polynomials, exponentials
      and power series
\endtitle
\author Jonathan Needleman\footnote{Research supported by the
            National Science Foundation through the Research
           Experiences for Undergraduates (REU) program at Cornell.
     Current address: Mathematics Department, Malott Hall, Cornell
         University, Ithaca, NY  14853.
                  \newline} \\
      {\rm Mathematics Department \\ Cornell University,
        Ithaca, NY 14853 \\ needlema\@math.cornell.edu } \\ \ \\
    Robert S. Strichartz\footnote{Research supported in part by the
         National Science Foundation, grant DMS-0140194.\newline} \\
      {\rm Mathematics Department, Malott Hall \\
       Cornell University, Ithaca, NY 14853 \\
      str\@math.cornell.edu} \\ \ \\
    Alexander Teplyaev \\ {\rm Mathematics Department \\
        University of Connecticut, Storrs, CT  06269\\
      teplyaev\@math.uconn.edu} \\ \ \\
    Po-Lam Yung\footnote{Research supported by the Mathematics Department
          of the Chinese University of Hong Kong, the Bankee Kwan Award
          for Mathematics Projects, and the Chung Chi Travelling Award
           for Mathematics.} \\
           {\rm Mathematics Department \\
         Chinese University of Hong Kong \\
       Shatin, Hong Kong, China \\ plyung\@math.cuhk.edu.hk}
\endauthor

\endtopmatter

\leftheadtext{J. Needleman, R. S. Strichartz, A. Teplyaev and P-L. Yung}

\rightheadtext{Calculus on the Sierpinski gasket}


\pageheight{7.5in}

\subheading{Abstract}
We study the analog of power series expansions on the Sierpinski gasket,
for analysis based on the Kigami Laplacian.  The analog of polynomials
are multiharmonic functions, which have previously been studied in
connection with Taylor approximations and splines.  Here the main
technical result is an estimate of the size of the monomials analogous to
$x^n/n!$. We propose a
definition of entire analytic functions as functions represented by power
series whose coefficients satisfy exponential growth conditions that are
stronger than what is required to guarantee uniform convergence.  We
present a characterization of these functions in terms of exponential
growth conditions on powers of the Laplacian of the function.  These
entire analytic functions enjoy properties, such as rearrangement and
unique determination by infinite jets, that one would expect.  However,
not all exponential functions (eigenfunctions of the Laplacian) are
entire analytic, and also many other natural candidates, such as the heat
kernel, do not belong to this class.  Nevertheless, we are able to use
spectral decimation to study exponentials, and in particular to
create exponentially decaying functions for negative eigenvalues.

\subheading{\S1.  Introduction}

\medskip
Ordinary calculus is such a remarkable subject because it combines both a
general conceptual framework and a detailed understanding of basic
functions.  For example, the theory of power series expansions hinges on
the elementary observation that the function $f_n(x) = x^n/n!$ on $[0,1]$
is bounded by $1/n!$.  (Stated this way, it seems almost a tautology, so
perhaps it is better to say that $f_n$ is the polynomial characterized
by the conditions $f_n^{(m)}(0) = \delta_{nm}$.)  Another example: among
all linear combinations of $\cosh x$ and $\sinh x$ there is one, $e^{-x}
= \cosh x -\sinh x$, that decays as $x \rightarrow \infty$; moreover its
rate of decay is the reciprocal of the growth rate of $\cosh x$ and $\sinh
x$.

The goal of this paper is to understand analogous facts about basic
functions on the Sierpinski gasket (SG), which should be regarded as the
simplest nontrivial example of a fractal supporting a theory of
differential calculus based on a Laplacian.  Standard references are the
books of Barlow [Ba] and Kigami [Ki2], and the expository paper [S2]. The
references to this paper, and the more extensive bibliography in [Ki2],
indicate an intensive development of the subject since Kigami's original
paper [Ki1] giving a direct analytic definition of the Laplacian on SG.

Recall that SG is the attractor of the iterated functions system (IFS)
consisting of three contractions in the plane $F_i(x) = \frac
12(x+q_i)$, $i = 0,1,2$ where $q_i$ are the vertices of an equilateral
triangle.  In other words $SG = \ds\bigcup^2_{i=0} F_i(SG)$, and we refer
to the sets $F_i(SG)$ as {\it cells of order} 1.  More generally, we
write $F_w = F_{w_1} \circ \cdots \circ F_{w_m}$ for a word $w =
(w_1,\ldots,w_m)$ of length $|w| = m$, each $w_j = 0, 1$ or $2$, and call
$F_w(SG)$ a {\it cell of level} $m$.  We regard SG as the limit of a
sequence of graphs $\Gamma_m$ (with vertices $V_m$ and edge relation $x
\sim_m y$) defined inductively as follows:  $\Gamma_0$ is the complete
graph on $V_0 = \{q_0,q_1,q_2\}$, and $V_m = \ds\bigcup^2_{i=0}
F_iV_{m-1}$ with $x \sim_m y$ if $x$ and $y$ belong to the same cell of
level $m$.  Then $V_* = \ds\bigcup^\infty_{m=1} V_m$, the set of all {\it
vertices}, the analog of the dyadic points in the unit interval, is dense
in SG.  We consider $V_0$ the set of boundary points of SG, and
$V_*\setminus V_0$ is the set of {\it junction points}.  Note that every
junction point in $V_m$ has exactly 4 neighbors in the graph $\Gamma_m$.
The graph Laplacian $\Delta_m$ on $\Gamma_m$ is defined by
$$
\Delta_m u(x) = \sum_{y\sim_m x} (u(y)-u(x)) \ \text{for} \
       x \in V_m \setminus V_0.    \tag1.1
$$
The Laplacian $\Delta$ on SG is defined as the renormalized limit
$$
\Delta u(x) = \lim_{m\rightarrow\infty} \frac 32 5^m\Delta_mu(x).
     \tag1.2
$$
More precisely, $u \in \dom \,\Delta$ and $\Delta u = f$ means $u$ and
$f$ are continuous functions and the limit on the right side of (1.2)
converges to $f$ uniformly on $V_*\setminus V_0$.  The Laplacian plays
the role of the second derivative on the unit interval (although it is
shown in [S] that it does not behave like a second order operator).  Thus
we will define a {\it polynomial} $P$ to be any solution of $\Delta^jP =
0$ for some $j$.  More precisely, if we let $\Cal H_j$ denote the space
of solutions of $\Delta^{j+1}u = 0$, then $\Cal H_j$ is a space of
dimension $3j+3$, and it has an ``easy'' basis $\{f_{nk}\}$ for $0 \le n
\le j$ and $k = 0,1,2$ characterized by
$$
\Delta^\ell f_{nk}(q_{k'}) = \delta_{\ell n} \delta_{kk'}.
     \tag1.3
$$
In [SU] a different basis was constructed in order to develop a theory of
splines.  Here we will consider yet another basis, implicitly used in [S3]
in conjunction with Taylor expansions, to define power series.

The Laplacian is basically an interior operator, as (1.2) is not defined
at the boundary (although $\Delta u = f$ makes sense at boundary points
by continuity).  There are also boundary derivatives.  The {\it normal
derivative}
$$
\partial_n u(q_j) = \lim_{m\rightarrow\infty} \Big(\frac 53\Big)^m
     (2u(q_j)-u(F^m_jq_{j+1}) - u(F^m_jq_{j-1}))
    \tag1.4
$$
(cyclic notation $q_{j+3} = q_j$) exists for every $u \in \dom\,\Delta$
and plays a crucial role in the theory, especially in the analog of the
Gauss-Green theorem:
$$
\int_{SG} (u\Delta v-v\Delta u)d\mu
      = \sum^2_{i=0} (u(q_i)\partial_nv(q_i) - \partial_nu(q_i)v(q_i)).
     \tag1.5
$$
Here $\mu$ is the natural probability measure that assigns weight
$3^{-m}$ to each cell of order $m$.  The normal derivative may be
localized to boundary points of any cell, and there is also a localized
version of (1.5).  At a junction point there are two different normal
derivatives with respect to the cells on either side.  For $u \in
\dom\,\Delta$ we have the {\it matching condition} that the two normal
derivatives sum to zero.  This leads to the {\it gluing property}: if $u$
and $f$ are continuous functions and $\Delta u = f$ on each cell of order
$m$ (meaning $\Delta (u\circ F_w) = 5^{-m}f\circ F_w$ for all words $w$
of length $m$), then $\Delta u = f$ on SG if and only if the matching
conditions hold at every junction point in $V_m$.

There are also {\it tangential derivatives}
$$
\partial_T u(q_j) = \lim_{m\rightarrow\infty} 5^m
       (u(F^m_0 q_{j+1}) - u(F^m_0q_{j-1}))
    \tag1.6
$$
that exist if $u \in \dom\,\Delta$, and may be localized to boundary
points of cells.  In this case there are no matching conditions for $u
\in \dom\,\Delta$.  However, we will show in Section 5 that there
are matching conditions involving infinite series of tangential and normal
derivatives valid for polynomials and analytic functions.  Tangential
derivatives were introduced in [S3].  Their true sigificance is still
somewhat elusive.  In this paper we will show that for polynomials and
analytic functions the sum of the tangential derivatives over the three
boundary points of any cell must vanish.  In [S3] and [T2] the idea of
creating a {\it gradient} of a function out of the normal and tangential
derivatives is discussed.  Here we will extend this to the idea of a {\it
jet}.  For simplicity we deal with a boundary point $q_\ell$,but the
definition can be localized to boundary points of any cell.

\medskip
\noindent
\underbar{Definition 1.1}:  For $u \in \dom\,\Delta^n$, the {\it $n$-jet}
of $u$ at $q_\ell$ is the $(3n+3)$-tuple of values
$(\Delta^ju(q_\ell),\partial_n\Delta^ju(q_\ell),\partial_T\Delta^ju(q_\ell))$
for $0 \le j \le n$. For $u \in \dom\,\Delta^\infty$,  the {\it jet} of
$u$ at $q_\ell$ is the infinite set of the same values for all $j \ge 0$.

Fix a boundary point $q_\ell$.  We define polynomials $P^{(\ell)}_{jk}$
by requiring that the $j$-jet at $q_\ell$ vanish except for one term,
$\Delta^jP^{(\ell)}_{j1}(q_\ell) = 1$,
$\partial_n\Delta^jP^{(\ell)}_{j2}(q_\ell) = 1$ and
$\partial_TP^{(\ell)}_{j3}(q_\ell) = 1$, respectively.  We refer to these
functions as {\it monomials}.  It is clear that the monomials
$P^{(\ell)}_{jk}$ for $0 \le j \le n$ form a basis of $\Cal H_n$.  It is
shown in [S] that they exhibit a prescribed decay rate in neighborhoods
of $q_\ell$, but the estimates established there were not uniform in
$j$.  The first goal of this paper is to obtain sharp estimates for
$\|P^{(\ell)}_{jk}\|_\infty$.   For $P^{(\ell)}_{j1}$ and
$P^{(\ell)}_{j3}$ we prove decay estimates faster than any exponential.
  For $P^{(\ell)}_{j2}$ the situation is different;  we prove an exponential
decay of order $\lambda_2^{-j}$ for the specific value $\lambda_2$ equal
to the second nonzero Neumann eigenvalue. This result is sharp.  In fact we
show that $(-\lambda_2)^jP_{j2}^{(\ell)}$ converges to a certain
$\lambda_2$--eigenfunction of $\Delta$.  This result has no analog in
ordinary calculus.

   We define a power
series about $q_\ell$ as an infinite linear combination of the monomials
$P^{(\ell)}_{jk}$ with coefficients $\{c_{jk}\}$.  We find growth
conditions on the coefficients to guarantee convergence.  We study the
rearrangement problem: given a convergent power seres about one boundary
point, does the function also have a convergent power series about the
other boundary points?  Surprisingly, we find that it is necessary to
assume a stronger growth restriction on the coefficients in order for
this to be the case, namely
$$
|c_{jk}| = O(R^j) \ \text{for some} \ R < \lambda_2.
   \tag1.7
$$
We end up defining an {\it entire analytic function} to be a function
represented by a power series with coefficients satisfying (1.7). We then
prove rearrangement is possible at all boundary points, and in fact local
power series expansions exist on all cells, with the estimate (1.7)
preserved (in fact the same $R$ value).  This choice of definition means
that there are some convergent power series that do not yield analytic
functions.  It also means that eigenfunctions of the Laplacian cannot be
entire analytic functions unless the eigenvalue satisfies $|\lambda| <
\lambda_2$.  On the other hand it is easy to see that there are
$\lambda_2$-eigenfunctions that cannot be represented by convergent power
series, so the definition seems to be close to best possible.  We then
are able to characterize the class of entire analytic functions in
$\dom\,\Delta^\infty$ by the growth conditions
$$
\|\Delta^ju\|_\infty = O(R^j) \ \text{for some} \ R < \lambda_2
    \tag1.8
$$
(one could also use $L^2$ norms).

Our definition of entire analytic function means that a basic principle
of unique analytic continuation holds.  If we have a function defined on
a cell and satisfying (1.8) there, it has a unique extension to an entire
analytic function on the whole space.  In fact its jet at any boundary
point of the cell satisfies (1.7), and uniquely determines the function.
This implies that a nonzero entire analytic function cannot vanish to
infinite order at any junction point.  We could also define local
analytic functions on a cell of order $m$ by relaxing the condition $R <
\lambda_2$ in (1.7) and (1.8) for $R < 5^m\lambda_2$.  One could hope to
have a notion of analytic continuation that would allow such local
analytic functions to extend to larger domains.  However, we have not been
able to find any interesting examples, so we will not pursue the matter
here.

It is easy to extend the notion of entire analytic function to infinite
blow--ups of SG ([S1], [T1]).  The simplest of these is
$$
SG_\infty = \bigcup^\infty_{n=1} F_0^{-n}(SG),  \tag1.9
$$
but more generally we could consider
$$
\bigcup^\infty_{n=1} F_{j_1}^{-1} F_{j_2}^{-1} \cdots F_{j_n}^{-1}(SG)
    \tag1.10
$$
for any choice of $j_1,j_2,j_3,\ldots$\ .  A function on SG satisfying
(1.8) for all $R > 0$ extends to an entire analytic function on any
blow--up (1.10).  It is not clear at present which, if any, of these
functions will come to play the role of special functions
(hypergeometric, Bessel functions, etc.) in real analysis.  On the other
hand it is very easy to construct many such functions simply by taking a
power series with bounded or sub--exponential growing coefficients.  The
negative results of [BST] mean that none of these spaces of analytic
functions is closed under multiplication, so this precludes using many
standard techniques for ordinary power series.

Although none of the eigenfunctions of the Laplacian are entire analytic
functions on the blow--ups, it is still important to understand their
global behavior.  In Section 6 we study this problem for the simplest
example $SG_\infty$ and negative eigenvalues.  It is easy enough to
define the analogs of the functions $\cosh\sqrt\lambda x$ and
$\sinh\sqrt\lambda x$.  In fact there are three, which we call
$C_\lambda(x)$, $S_\lambda(x)$ and $Q_\lambda(x)$, characterized among
$(-\lambda)$--eigenfunctions by their $0$--jet at $q_0$, or equivalently
by power series involving just $P_{j1}^{(0)}$, $P_{j2}^{(0)}$, or
$P_{j3}^{(0)}$ terms, respectively.  The power series for $C_\lambda(x)$
and $Q_\lambda(x)$ converge on all of $SG_\infty$, while the power series
for $S_\lambda(x)$ is only convergent on a neighborhood of $q_0$
(depending on $\lambda$).  Fortunately, there is another method available
to study these eigenfunctions, called {\it spectral decimation} ([FS],
[DSV], [T1]).  Using this method we are able to show that they exhibit an
exponential growth as $x \rightarrow \infty$ (or as $\lambda \rightarrow
\infty$), and there is one linear combination, $E_\lambda(x) =
C_\lambda(x)-S_\lambda(x)$ for the appropriate normalization, that decays
as $x \rightarrow \infty$ at the reciprocal rate.  Thus $E_\lambda(x)$ is
the analog of $e^{-\sqrt\lambda x}$.  It is not clear if there is any
analog of $e^{\sqrt\lambda x}$.

Although we do not use power series in our study of properties of
eigenfunctions, we can turn the tables and use facts about eigenfunctions
to obtain information about power series.  In particular, we are able to
construct specific power series that are divergent, or power series that
are convergent but not rearrangeable.  We can also give an explanation
for why the recursion relations for the size of monomials are unstable.

It is interesting to speculate on possible future extensions and
developments of our results.  It is important to understand all
eigenfunctions, including those with positive eigenvalues, on all
blow--ups (1.10).  There should be some sort of Liouville--type theorem
precluding nonconstant bounded entire analytic functions on blow--ups
without boundary.

What is the behavior of an entire analytic function in a neighborhood of
a generic point?  Is there any notion of power series there?  Are there
interesting examples of local analytic functions with a natural domain
that is not just a single cell?  Is there a meaningful notion of analytic
functions on fractafolds based on SG [S4]?

We have seen that there is no restriction on the jet of an analytic
function other than the growth condition (1.7).  For the larger class
$\dom \,\Delta^\infty$, is there an analog of Borel's theorem that an
arbitrary jet may be specified at one (or all three) boundary points?

In [OSY], the structure of level sets of harmonic functions on SG was
elucidated, with the remark that certain eigenfunctions of the Laplacian
have level sets of an entirely different nature.  It is clear now that
these eigenfunctions are not analytic, so it is reasonable to ask if
anything interesting can be said about level sets of entire analytic
functions.  Another remark from that paper is that harmonic functions
enjoy a principle called ``geography is destiny.''  Roughly speaking,
this says that the restriction to a small cell of a harmonic function is
essentially dictated (up to two parameters) by the location of the cell,
rather than the specific harmonic function, in a certain generic sense.
This holds because restrictions of harmonic functions are governed by
long products of matrices, so the theory of products of random matrices
makes generic predictions.  For analytic functions, there is a similar
description of the transformation of jets, except that the matrices are
now infinite.  So if we go to a small cell, while all jets satisfying
(1.7) are possible, some may be very unlikely for a generic analytic
function.  Is there some way to make this precise?

A sequel to this paper, [BSSY], will discuss functions with point
singularities, exponential functions on general blow--ups, and estimates
for normal derivatives of Dirichlet eigenfunctions and heat kernels.

\newpage
\subheading{\S2.  Polynomials}

\medskip
The space $\Cal H_j$ of $(j+1)$-harmonic functions (solutions of
$\Delta^{j+1}u = 0$) has dimension $3(j+1)$ and plays the role of the
space of polynomials of degree at most $2j+1$ on the unit interval.
Several different bases for $\Cal H_j$ are known.  In [SU], in order to
develop a theory of spline spaces,  bases based on the
behavior at all three boundary points were used. In this section we
will discuss properties of yet another basis, based on the behavior at a
single boundary point, that is more suited to the work on power series to
follow.  The polynomials in this basis are analogous to the monomials
$x^n/n!$ on the unit interval.  These functions were introduced in [S3],
but not much was done there to describe their behavior.

\medskip
\noindent
\underbar{Definition 2.1}:  Fix a boundary point $q_\ell$.  The
{\it monomials} $P_{jk}^{(\ell)}$ for $k = 1,2,3$ and $j =
0,1,2,\ldots$ are defined to be the functions in $\Cal H_j$ satisfying
$$
\Delta^m P^{(\ell)}_{jk}(q_\ell) =  \delta_{mj}\delta_{k1}
     \tag2.1
$$
$$
\partial_n \Delta^m P^{(\ell)}_{jk} (q_\ell) = \delta_{mj}\delta_{k2}
    \tag2.2
$$
$$
\partial_T \Delta^m P^{(\ell)}_{jk}(q_\ell) = \delta_{mj}\delta_{k3}.
     \tag2.3
$$
When $\ell = 0$ we will sometimes delete the upper exponent and just
write $P_{jk}$.

Note that we only need to consider $m \le j$ in (2.1-3), since
$\Delta^mP^{(\ell)}_{jk}$ vanishes identically otherwise.  Thus there are
$3(j+1)$ conditions in all, and it follows from [S3] that there is a
unique solution, and the monomials $P^{(\ell)}_{jk}$ for fixed $\ell$ and
all $j \le j_1$ form a basis for $\Cal H_{j_1}$.  We have the self-similar
identities
$$
P^{(\ell)}_{j1}(F^m_\ell x) = 5^{-jm} P_{j1}^{(\ell)}(x)
    \tag2.4
$$
$$
P^{(\ell)}_{j2}(F^m_\ell x) = \Big(\frac 35\Big)^m 5^{-jm}
         P^{(\ell)}_{j2}(x)
      \tag2.5
$$
$$
P^{(\ell)}_{j3} (F^m_\ell x) = 5^{-(j+1)m} P_{j3}^{(\ell)}(x)
       \tag2.6
$$
that describe the decay rate of these functions as $x \rightarrow q_\ell$
(of course $P^{(\ell)}_{01} \equiv 1$).  It is easy to see that
$P^{(\ell)}_{j1}$ and $P^{(\ell)}_{j2}$ are symmetric while
$P^{(\ell)}_{j3}$ is skew-symmetric under the reflection that fixes
$q_\ell$ and permutes the other two boundary points.  It is easy to
compute the values of monomials to any desired precision.  Figure 2.1
shows the graphs of some of them.  Since we may obtain $P^{(\ell)}_{jk}$
from $P^{(0)}_{jk}$ by simply rotating the variable $x$, we will
restrict our discussion to $\ell = 0$ from now on.

\newpage
\

\bigskip

\centerline{
\BoxedEPSF{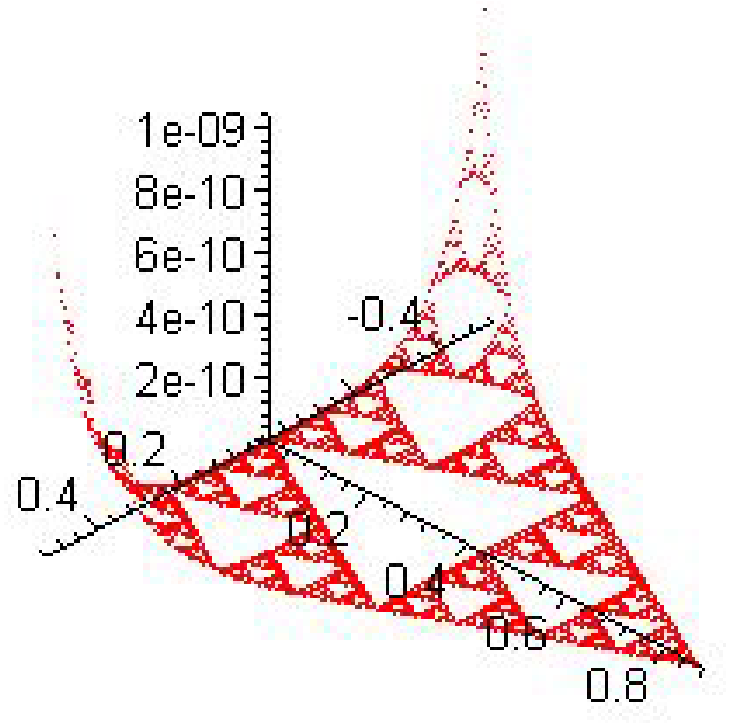 scaled 700} \hskip.4in \BoxedEPSF{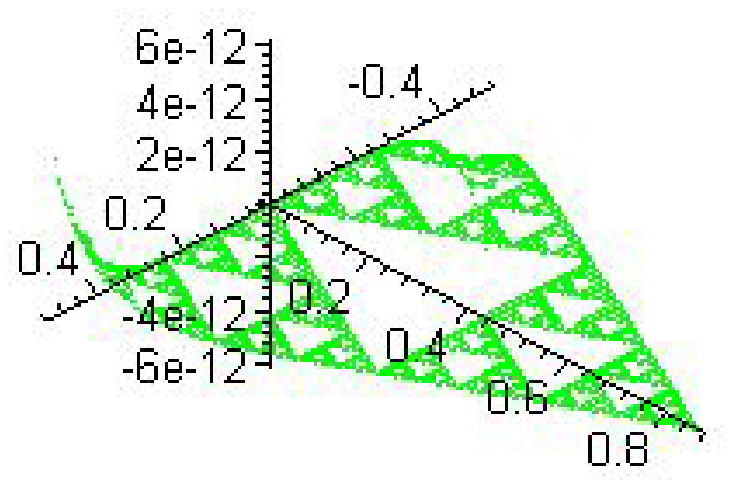
scaled 700}}

\bigskip
\centerline{
\BoxedEPSF{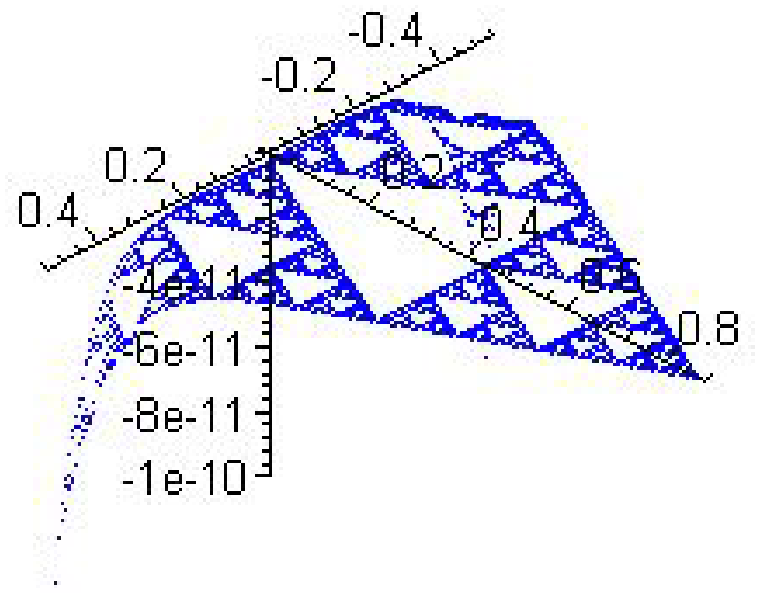 scaled 700} \hskip.4in \BoxedEPSF{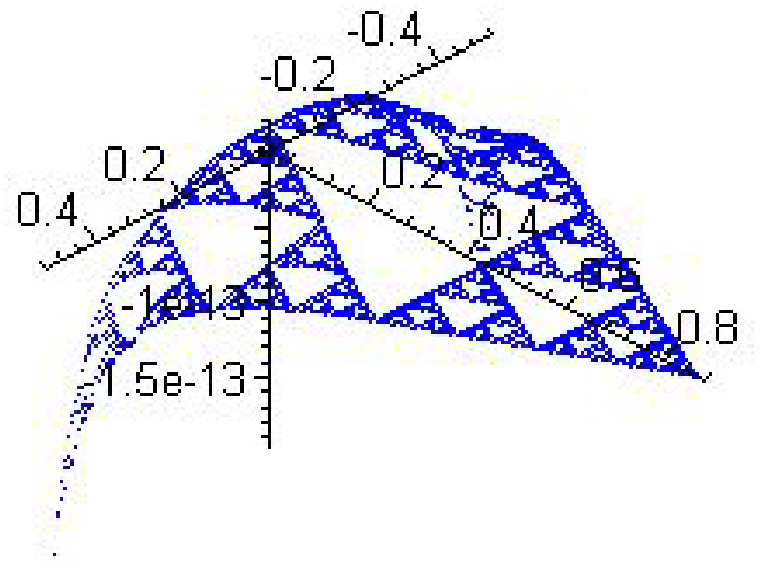
scaled 700}}

\bigskip
\centerline{
\BoxedEPSF{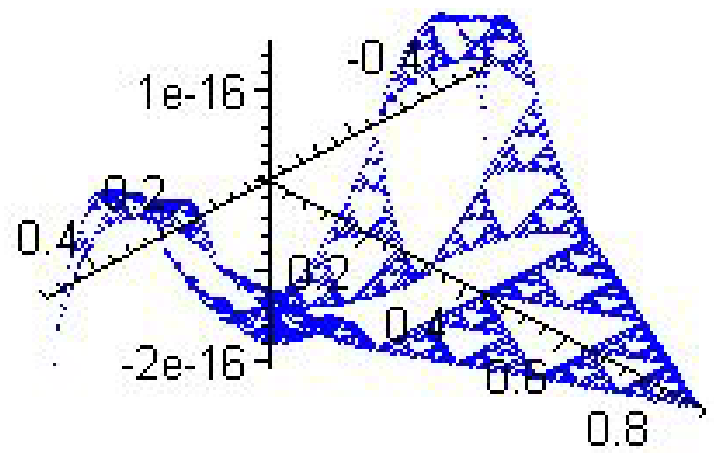 scaled 700} \hskip.4in \BoxedEPSF{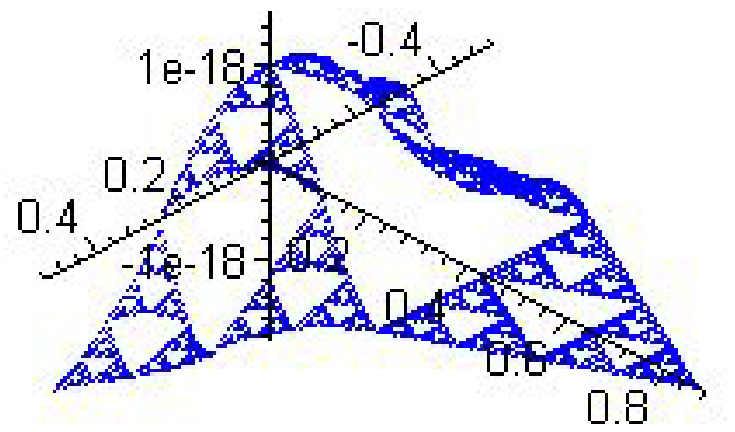
scaled 700}}

\bigskip
\noindent
\underbar{Figure 2.1}:  The graphs of $P_{jk}$ for some typical values.
The graphs of $P_{j1}$ are all qualitatively similar for $j \ge 1$, so we
show only $P_{51}$.  Similarly for $P_{j3}$.  The nature of the graphs of
$P_{j2}$ changes drastically around $j = 5,6,7,8$, so we display all of
these.  The graphs of $P_{j2}$ for $j \ge 8$ are qualitatively similar to
$P_{82}$.


\newpage

It is clear from the definition that powers of the Laplacian send
monomials to monomials, simply reducing the $j$ index:
$$
\Delta^m P_{jk} = P_{(j-m)k}.
    \tag2.7
$$
We could use this property to give an inductive definition.  When $j = 0$
the monomials are explicit harmonic functions, $P_{01} \equiv 1$,
$P_{02}$ has boundary values $P_{02}(q_0) = 0$, $P_{02}(q_1) =
P_{02}(q_2) = -1/2$ and $P_{03}$ has boundary values $P_{03}(q_0) = 0$,
$P_{03}(q_1) =~\!-P_{03}(q_2) = 1/2$.  Then $P_{jk}$ for $j > 0$ is the
unique solution of $\Delta P_{jk} = P_{(j-1)k}$ with vanishing initial
conditions
$$
P_{jk}(q_0) = 0, \ \partial_n P_{jk}(q_0) = 0, \
        \partial_T P_{jk}(q_0) = 0.
$$
In [KSS] it is shown that $P_{jk}$ may then be written as an integral
operator (with explicit kernel) applied to $P_{(j-1)k}$.  However, the
kernel is quite singular, so we have not been able to extract any useful
information out of this representation.

There are three main goals in this section: 1) to obtain sharp estimates
for the size of the monomials, 2) to understand how to express monomials
for one choice of $\ell$ in terms of monomials for another choice of
$\ell$, 3) to obtain certain universal identities that hold for all
monomials.  In pursuit of these goals we introduce some terminology.

\medskip
\noindent
\underbar{Definition 2.2}:  For $j \ge 0$ let
$$
\cases
   &\alpha_j = P_{j1}(q_1), \ \beta_j = P_{j2}(q_1), \
                     \gamma_j = P_{j3}(q_1) \\
  &n_j = \partial_n P_{j1}(q_1), \ t_j = \partial_TP_{j2}(q_1).
\endcases     \tag2.8
$$

Note that by symmetry we have $P_{j1}(q_2) = \alpha_j$, $P_{j2}(q_2) =
\beta_j$ and $P_{j3}(q_2) = -\gamma_j$, so that all values of monomials
at boundary points are expressible in terms of $\alpha$'s, $\beta$'s and
$\gamma$'s.  Soon we will see that the $n$'s, $t$'s and $\alpha$'s
suffice to express all normal and tangential derivatives of monomials at
boundary points.

\proclaim{Theorem 2.3}  The following recursion relations hold:
$$
\alpha_j = \frac{4}{5^j-5} \sum^{j-1}_{\ell=1} \alpha_{j-\ell}
              \alpha_\ell \quad \text{for} \quad j \ge 2
      \tag2.9
$$
$$
\gamma_j = \frac{4}{5^{j+1}-5} \sum^{j-1}_{\ell=0} \alpha_{j-\ell}
           \gamma_\ell \quad \text{for} \quad j \ge 1
    \tag2.10
$$
$$
\beta_j = \frac{1}{5^j-1} \sum^{j-1}_{\ell=0} \Big(\frac 25 5^{j-\ell}
                \alpha_{j-\ell}\beta_\ell - \frac 23 \alpha_{j-\ell}
            5^\ell \beta_\ell + \frac 45 \alpha_{j-\ell} \beta_\ell\Big)
         \quad \text{for} \quad j \ge 1,
     \tag2.11
$$
with initial data $\alpha_0 = 1$, $\alpha_1 = 1/6$, $\beta_0 = -1/2$,
$\gamma_0 = 1/2$.  In particular,
$$
\gamma_j = 3\alpha_{j+1}.
     \tag2.12
$$
\endproclaim

\underbar{Proof}:  It is convenient to work in matrix notation, with all
matrices being infinite semi-circulant.  For example, the matrix $\alpha
= \{\alpha_{ij}\}_{i,j=0,1,2,\ldots}$ has $\alpha_{ij} = \alpha_{i-j}$
for $i \ge j$ and $\alpha_{ij} = 0$ for $i < j$.  We consider two linear
operators on such matrices, the shift $\sigma$ and the dilation $\tau$,
given by
$$
\sigma \pmatrix  d_0  &0  &\cdots \\ d_1  &d_0  & 0 \\
             d_2  &d_1  &d_0  &0 \\      \vdots
        \endpmatrix
= \pmatrix  d_1  &0  &\cdots \\ d_2  &d_1  & 0 \\
             d_3  &d_2  &d_1  &0 \\      \vdots
        \endpmatrix
$$
$$
\tau \pmatrix  d_0  &0  &\cdots \\ d_1  &d_0  & 0  &\cdots\\
             d_2  &d_1  &d_0  &0 &\cdots \\      \vdots
        \endpmatrix
=   \pmatrix d_0 &0   &\cdots \\  5d_1  &d_0  &0 \\
        5^2d_2  &5d_1  &d_0  &0 \\  \vdots
     \endpmatrix .
$$

Let $\{f_{j1},f_{j2},f_{j3}\}_{j=0}^\infty$ be the easy basis defined by
(1.3). As in [SU] we let
$$
\split
     a_{l-1}&=\partial_n f_{lk}(q_k) \\
     b_{l-1}&=\partial_n f_{lk}(q_n) \qquad n \ne k
\endsplit
$$
for $l=0,1,2,\dots$. Then the Gauss-Green formula says for $l\geq 0$
$$
\split
     a_l &=\partial_n f_{(l+1)1}(q_1)\\
         &=\sum_{n=1}^3\left(f_{01}(q_n)\partial_n
         f_{(l+1)1}(q_n)-f_{(l+1)1}(q_n) \partial_n f_{01}(q_n)\right)\\
         &=\int_{SG} (f_{01} \Delta f_{(l+1)1}- f_{(l+1)1} \Delta f_{01})
         d\mu \\
         &=\int_{SG} f_{01}f_{l1} d\mu
\endsplit
$$
and
$$
\split
     b_l &=\partial_n f_{(l+1)1}(q_2)\\
         &=\sum_{n=1}^3\left(f_{02}(q_n)\partial_n
         f_{(l+1)1}(q_n)-f_{(l+1)1}(q_n) \partial_n f_{02}(q_n)\right)\\
         &=\int_{SG} (f_{02} \Delta f_{(l+1)1}- f_{(l+1)1} \Delta f_{02})
         d\mu \\
         &=\int_{SG} f_{02}f_{l1} d\mu.
\endsplit
$$
This shows that our definition is consistent with [SU].  It is easy to see
that $a_{-1}=2, b_{-1}=1$.

We note here some typos from [SU]:

\medskip
(i) in (5.4) the coefficient $\frac{47}{45}$ should be
$\frac{47}{75}$;

\medskip
(ii) in the first line of (5.7) the coefficients 2 of
$a_{j-1-\ell}$ and
$b_{j-1-\ell}$ should be deleted.

Now let $p_j$, $q_j$ be defined by
$$
\split
&p_j = 5^jf_{jk}(F_iq_k) \quad i \neq k \\
&q_j = 5^jf_{jk}(F_iq_\ell) \quad \text{for} \ i, j, \ell \
        \text{distinct}.
\endsplit
$$
(Note that we are using the same symbol $q_j$ for two different things,
but it should be clear from context which is which.)

Then (5.7) of [SU] rearranged says
$$
\split
     \sum_{l=0}^j (a_{j-l-1} +b_{j-l-1}) (2p_l+q_l) + b_{j-1}&=0    \\
     \sum_{l=0}^j (2a_{j-l-1}-b_{j-l-1}) (p_l-q_l) + b_{j-1}&=0
\endsplit
$$

If we set
$$
A  =  \pmatrix
    a_{-1} & 0  \\
   a_0 & a_{-1} & 0 \\
   a_1 & a_0 & a_{-1} & 0  \\
   a_2 & a_1 & a_0 & a_{-1} & \ddots \\
   \vdots &  &  & \ddots  & \ddots
    \endpmatrix \qquad
B = \pmatrix
   b_{-1} & 0  \\
   b_0 & b_{-1} & 0 \\
   b_1 & b_0 & b_{-1} & 0  \\
   b_2 & b_1 & b_0 & b_{-1} & \ddots \\
   \vdots &  &  & \ddots  & \ddots
\endpmatrix
$$
$$
P =
   \pmatrix
   p_0 & 0  \\
   p_1 & p_0  & 0 \\
   p_2 & p_1 & p_0 & 0  \\
   p_3 & p_2 & p_1 & p_0 & \ddots \\
   \vdots &  &  & \ddots  & \ddots
\endpmatrix    \qquad
Q = \pmatrix
   q_0 & 0  \\
   q_1 & q_0  & 0 \\
   q_2 & q_1 & q_0 & 0  \\
   q_3 & q_2 & q_1 & q_0 & \ddots \\
   \vdots &  &  & \ddots  & \ddots
\endpmatrix .
$$
Then in matrix notation this becomes

$$
(A+B)(2P+Q)+B = 0,  \qquad
(2A-B)(P-Q)+B = 0. \tag2.13
$$

Now for $j \geq 0$,
$$
\cases
     P_{j1}= f_{j0} + \ds\sum_{l=0}^j \alpha_{j-l} (f_{l1}+f_{l2}) \\
     P_{j2}=\ds\sum_{l=0}^j \beta_{j-l} (f_{l1}+f_{l2}),
\endcases  \tag2.14
$$
so taking normal derivatives at $q_0$, we have
$$
\split
     a_{j-1} + 2 \sum_{l=0}^j \alpha_{j-l}b_{l-1}&=\partial_n
               P_{j1}(q_0)=0\\
     2\sum_{l=0}^j \beta_{j-l} b_{l-1} &= \partial_n P_{j2} (q_0)
   = \cases 1, & \text{if}\  $j=0$; \\
     0, & \text{otherwise}. \endcases
\endsplit
$$
In matrix notation this is
$$  2\alpha B +A=0, \quad 2\beta B = I,
$$
i.e.
$$ A=-\alpha \beta^{-1}, \qquad
     B=\frac{1}{2}\beta^{-1}  \tag2.15
$$

Substituting (2.15) into (2.13), we get
$$
\split
     2P+Q &= -(A+B)^{-1} B = -[-\frac{1}{2} \beta^{-1} (2\alpha-I)]^{-1}
             [\frac{1}{2} \beta^{-1}] = (2\alpha-I)^{-1} \\
     P-Q &= -(2A-B)^{-1} B = -[-\frac{1}{2} \beta^{-1} (4\alpha+I)]^{-1}
             [\frac{1}{2} \beta^{-1}] = (4\alpha+I)^{-1}
\endsplit
$$
so
$$
     (2\alpha-I)(2P+Q)=I=(4\alpha+I)(P-Q).
$$
Expanding we get
$$
     4\alpha P+2\alpha Q-2P-Q=4\alpha P-4\alpha Q+P-Q,
$$
i.e.
$$
   P=2\alpha Q, \ \text{and} \quad
     Q=(4\alpha +I)^{-1}(2\alpha-I)^{-1}.  \tag2.16
$$

Now evaluate (2.14) at $F_0q_1$, noting
that
$$
\split
     P_{j1}(F_0q_1)&=5^{-j}P_{j1}(q_1)=5^{-j}\alpha_j\\
     P_{j2}(F_0q_1)&=\frac{3}{5}5^{-j}P_{j1}(q_1)=\frac{3}{5}
     5^{-j} \beta_j
\endsplit
$$
by (2.4), (2.5) and
$$
\aligned
     f_{l0}(F_0q_1)&=f_{l1}(F_0q_1)=5^{-l}p_l \\
     f_{l2}(F_0q_1)&=5^{-l}q_l ,
\endaligned    \tag2.17
$$
by the definitions of $p_l$'s and $q_l$'s.  The result is
$$
\split
     &5^{-j}\alpha_j =5^{-j}p_j+\sum_{l=0}^j \alpha_{j-l}
     (5^{-l}p_l+5^{-l}q_l)\\
     &\frac{3}{5}5^{-j}\beta_j=\sum_{l=0}^j \beta_{j-l}
     (5^{-l}p_l+5^{-l}q_l)
\endsplit
$$
so
$$
\split
     &\alpha_j = p_j+\sum_{l=0}^j 5^{j-l}\alpha_{j-l} (p_l+q_l) \
         \text{and}\\
     &\frac{3}{5}\beta_j = \sum_{l=0}^j 5^{j-l} \beta_{j-l} (p_l+q_l).
\endsplit
$$
In matrix notation these read
$$
\split
     &\alpha = P+ \tau(\alpha) (P+Q) \ \text{and} \\
     &\frac{3}{5}\beta = \tau(\beta)(P+Q).
\endsplit
$$
 From (2.14) we see that
$$
\split
     &\alpha=[2\alpha+\tau(\alpha)(2\alpha+I)]Q \ \text{and}\\
     &\frac{3}{5} \beta = \tau(\beta)(2\alpha+I)Q
\endsplit
$$
hence
$$
\split
     &\tau(\alpha) = 4\alpha^2 - 3\alpha \ \text{and}\\
     &\frac{3}{5}\beta(2\alpha-I)(4\alpha+I)=\tau(\beta)(2\alpha+I),
\endsplit
$$
from which (2.9) and (2.11) follow.

Finally
$$
     P_{j3} = \sum_{l=0}^j \gamma_{j-l}(f_{l1}-f_{l2})
$$
$$
     P_{j3}(F_0q_1) = 5^{-(j+1)} P_{j3}(q_1) = 5^{-(j+1)} \gamma_j
$$
and so by (2.17) we have
$$
     5^{-(j+1)}\gamma_j = \sum_{l=0}^j \gamma_{j-l}
     (5^{-l}p_l-5^{-l}q_l),
$$
i.e.
$$
     \frac{1}{5}\gamma_j = \sum_{l=0}^j 5^{j-l}\gamma_{j-l} (p_l - q_l),
$$
or in matrix notation
$$
     \frac{1}{5}\gamma = \tau(\gamma) (P-Q).
$$
Thus $\tau(\gamma) = \frac 15(4\alpha+I)\gamma$  from which
(2.10) follows.

The values of $\alpha_0$, $\beta_0$ and $\gamma_0$ are easy to check.
Then (2.12) follows from (2.9) and (2.10) since $\alpha_j$ and
$\alpha_{j-1}$ satifsy the same recursion relation.
\hfill Q.E.D.

\bigskip
\proclaim{Theorem 2.4}  For all $j \ge 0$ we have
$$
P^{(0)}_{j3}(x) + P^{(1)}_{j3}(x) + P^{(2)}_{j3}(x) = 0
    \tag2.18
$$
and
$$
P^{(0)}_{j3}(x) = 3(P^{(2)}_{(j+1)1}(x) - P^{(1)}_{(j+1)1}(x)).
    \tag2.19
$$
\endproclaim

\underbar{Proof}:  We prove (2.18) by induction.  For $j = 0$ the left
side is a harmonic function that vanishes on the boundary (because of the
skew-symmetric of each term). Such a function must be zero.  For
the induction step, assume it is true for $j-1$.  Then
$$
\Delta(P^{(0)}_{j3} + P^{(1)}_{j3} + P^{(2)}_{j3})
      = P^{(0)}_{(j-1)3} + P^{(1)}_{(j-1)3} + P^{(2)}_{j3} = 0
$$
by the induction hypothesis.  Once again the left side is a harmonic
function, and it vanishes on the boundary by skew symmetry.

To prove (2.19) we use
$$
P^{(0)}_{j3} = \sum^j_{\ell=0} \gamma_{j-\ell}
          (f_{\ell 1}-f_{\ell 2}).  \tag2.20
$$
On the other hand, we have
$$
\split
&P^{(2)}_{(j+1)1} = f_{(j+1)2} + \sum^{j+1}_{\ell=0}
              \alpha_{j-\ell+1}(f_{\ell 0} + f_{\ell 1}) \\
   &P^{(1)}_{(j+1)1} = f_{(j+1)1} + \sum^{j+1}_{\ell=0}
          \alpha_{j-\ell+1} (f_{\ell0} + f_{\ell 2})
\endsplit
$$
so that
$$
\split
P^{(2)}_{(j+1)1} - P^{(1)}_{(j+1)1}
     &= f_{(j+1)2} - f_{(j+1)1} + \sum^{j+1}_{\ell=0} \alpha_{j-\ell+1}
           (f_{\ell 1} - f_{\ell 2}) \\
   &= \sum^j_{\ell=0} \alpha_{j-\ell+1} (f_{\ell 1} - f_{\ell 2})
\endsplit
$$
since $\alpha_0 = 1$.  The result follows from (2.12).  \hfill  Q.E.D.

\bigskip
The dihedral-3 symmetry group $D_3$ of SG consists of reflections
$\rho_0$, $\rho_1$, $\rho_2$, where $\rho_j$ preserves $q_j$ and permutes
the other two boundary points, and rotations $I$, $R_1$, $R_2 = (R_1)^2$
where $R_1 q_j = q_{j+1}$ (cyclic notation).

\proclaim{Theorem 2.5}  Any polynomial $P$ satisfies the identity
$$
P(x) + P(R_1x) + P(R_2x) = P(\rho_0x) + P(\rho_1x) + P(\rho_2x),
    \tag2.21
$$
and more generally the local versions
$$
P(x_0) + P(x_1) + P(x_2) = P(y_1) + P(y_2) + P(y_3)
    \tag2.22
$$
for any sextuplet of points such that
$$
\cases
   &x_0 = F_wx, \ x_1 = F_wR_1x, \ x_2 = F_w R_2x, \\
  &y_0 = R_w\rho_0x, \ y_1 = F_w \rho_1 x, \ y_2 F_w \rho_2x
\endcases    \tag2.23
$$
for some $x \in SG$ and some word $w$.
\endproclaim

\underbar{Proof}:  The local version follows from (2.21) because $P \circ
F_w$ is also a polynomial.  To prove (2.21) it suffices to show it holds
for all monomials.  Now we claim (2.21) is trivially true for any
function that is symmetric with respect to one of the reflections
$\rho_j$.  Say $P(x) = P(\rho_0x)$ for all $x$.  Then $P(R_1x) =
P(\rho_1x)$ and $P(R_2x) = P(\rho_2x)$ because $\rho_0R_1 = \rho_1$ and
$\rho_0R_2 = \rho_2$.  In particular, (2.21) holds for all
$P^{(\ell)}_{j1}$ and $P^{(\ell)}_{j2}$.  It follows from (2.19) that
it also holds for $P^{(\ell)}_{j3}$.   \hfill  Q.E.D.

\bigskip
The same result holds for uniform limits of polynomials; in particular,
the convergent power series discussed in the next section.  Note that
Kigami [Ki2] Theorem 4.3.6 has characterized the space of $L^2$ limits of
polynomials by the condition of orthogonality to all joint Dirichlet and
Neumann eigenfunctions.  It is not hard to see that (2.22) implies the
orthogonality to some of these eigenfunctions (those of the
$\lambda^{(5)}$-type in [DSV]), but not others.  On the other hand, it is
not clear how these orthogonality conditions imply (2.22).

\bigskip
\proclaim{Corollary 2.6} Any polynomial $P$ satisfies
$$
\partial_T P(q_0) + \partial_T P(q_1) + \partial_T P(q_2) = 0,
     \tag2.24
$$
and more generally the sum of tangential derivtives at the boundary
points of any cell must vanish.
\endproclaim

\underbar{Proof}:  Taking $x = F^m_0q_1$ in (2.21), we find
$$
(P(F^m_0q_1) - P(F^m_0q_2)) + (P(F^m_1q_2) - P(F^m_1q_0))
     + (P(F^m_2q_0) - P(F^m_2q_1)) = 0
    \tag2.25
$$
because $R_1F^m_0q_1 = F^m_1q_2$, $R_2F^m_0q_1 = F^m_2q_0$,
$\rho_0F_0^mq_1 = F^m_0q_2$, $\rho_1F^m_0q_1 = F^m_2q_1$, $\rho_2F^m_0q_1
= F^m_1q_0$.  Multiplying (2.25) by $5^m$ and taking the limit as $m
\rightarrow \infty$ yields (2.24).  The local form follows as before.
\hfill  Q.E.D.

\bigskip
\underbar{Remark}:  As we observed in the proof of Theorem 2.5, any
polynomial may be written as a sum of three polynomials, each symmetric
with respect to one of the reflections $\rho_j$, $P = P^{(0)} + P^{(1)} +
P^{(2)}$.  It is easy to see that one way to do this explicitly is to take
$$
P^{(j)}(x) = \frac 13 (P(x) + P(\rho_jx)) - \frac 19 (P(\rho_0x)
           + P(\rho_1x) + P(\rho_2x)).
    \tag2.26
$$

We consider next estimates for the size of $\alpha_j$, $\beta_j$,
$\gamma_j$.  We show that $\alpha_j$ has rapid decay, which we believe is
fairly sharp.  This gives the same decay rate for $\gamma_j$.

\proclaim{Theorem 2.7}  There exists a constant $c$ such that
$$
0 < \alpha_j < c(j!)^{-\log 5/\log 2} \quad \text{for all} \quad j.
    \tag2.27
$$
\endproclaim

\underbar{Proof}:  It is clear from (2.9) and the initial conditions that
the $\alpha_j$ are positive.  Let $\wt\alpha_j = (j!)^{\log 5/\log
2}\alpha_j$.  We need to show that the $\wt\alpha_j$ are bounded, which
we do by induction.  If $\wt\alpha_\ell \le c$ for $\ell \le j$, then
(2.9) implies
$$
\wt\alpha_j \le c^25^{1-j} \sum^{j-1}_{\ell=1}
         \binom{j}{\ell}^{\log 5/\log 2}.
$$
It is well known that
$$
\sum^j_{\ell=0} \binom{j}{\ell}^2 = \binom{2j}{j},
$$
so by Stirling's formula and routine arguments we have
$$
\sum^{j-1}_{\ell=1} \binom{j}{\ell}^{\log 5/\log 2}
     \le M5^j(j)^{-1/2}
$$
for all $j \ge 2$ for a small constant $M$, so $\wt\alpha_j \le
c^25M(j)^{-1/2}$.  It is easy to choose $c$ and $j_0$ so that
$\wt\alpha_\ell \le c$ for $\ell < j_0$ and $c \le (j_0)^{1/2}/5M$.
\hfill  Q.E.D.

\bigskip
Table 2.1 presents numerical computations of $\alpha_j$ and $\beta_j$.

\newpage

$$
\matrix
j     &\alpha_j   &\beta_j   &(-\lambda_2)^j\beta_j
       &8^j(j!)^{\frac{\log(5)}{\log(2)}}\alpha_j\\
0  &1.   &-.5000000000 &-.5000000000 &1. \\
1  &.1666666667 &-.04444444444 &6.025427867 &1.333333333 \\
2  &.005555555556 &-.001008230453 &-18.53107571 &1.777777777 \\
3  &.00006172839506 &-.8554950809 \   10^{-5} &21.31713060 &2.025658338 \\
4  &.3318730917 \ 10^{-6} &-.3853047646 \ 10^{-7} &-13.01625411
               &2.178127244 \\
5  &.1021147975\  10^{-8} &-.9848282711
               \ 10^{-10} &4.510374011 &2.250339083 \\
6  &.2007235906 \ 10^{-11} &-.1933836698\  10^{-12}
              &-1.200721414 &2.268082964
\\ 7  &.2713115918\  10^{-14} &-.7720311754 \ 10^{-16} &.06498718216
                        &2.248411184
\\ 8  &.2656437390\  10^{-17} &-.1187366658\  10^{-17} &-.1355027558
               &2.201440598
\\ 9  &.1959165201 10^{-20} &.7232200062\  10^{-20} &-.1118933095
                &2.134277683
\\ 10 &.1122370097 \ 10^{-23} &-.5436238235 \ 10^{-22} &-.1140256558
                 &2.052740417
\\ 11 &.5120236416 \ 10^{-27} &.4004514705 \ 10^{-24} &-.1138739539
                  &1.961629028
\\ 12 &.1898528071 \  10^{-30} &-.2954013973 \ 10^{-26} &-.1138826233
                    &1.864726441
\\ 13 &.5820142006 \ 10^{-34} &.2178916451\  10^{-28} &-.1138822148
                  &1.764891613
\\ 14 &.1496625756\  10^{-37} &-.1607201123 \ 10^{-30} &-.1138822304
                   &1.664234594
\\ 15 &.3268360869 \ 10^{-41} &.1185495242 \ 10^{-32} &-.1138822298
                     &1.564302197
\\ 16 &.6126918156\  10^{-45} &-.8744387717 \ 10^{-35} &-.1138822298
                 &1.466232140
\\ 17 &.9952451630 \ 10^{-49} &.6449989323\  10^{-37} &-.1138822298
                &1.370864839
\\ 18 &.1412543698 \ 10^{-52} &-.4757607235\  10^{-39} &-.1138822298
                &1.278818576
\\ 19 &.1764707126 \ 10^{-56} &.3509281252 \ 10^{-41} &-.1138822298
                 &1.190538877
\\ 20 &.1953558627 \ 10^{-60} &-.2588497599 \ 10^{-43} &-.1138822298
                    &1.106332006
\\
\endmatrix
$$

\centerline{Table 2.1.}
\bigskip
\noindent
It appears that $8^j(j!)^{\log 5/\log 2} \alpha_j$ remains bounded (8 is by
no means the best constant, and perhaps it could be replaced by an
arbitrary positive number).   It also appears that $(-\lambda_2)^j\beta_j$
converges to the constant $-.1138822298$, where $\lambda_2 =
135.572126995788\ldots$ is the  second nonzero Neumann eigenvalue.  It is
easy to see that $\lambda_2$ is the largest value for which such an
estimate could hold, because
$$
\sum^\infty_{j=0} \beta_j(-\lambda_2)^j \ \text{diverges}.
$$
Indeed, if we did not have divergence then
$$
\sum^\infty_{j=0} (-\lambda_2)^j P_{j2}(x)
$$
would be a solution to the eigenvalue equation $-\Delta u = \lambda_2u$
satisfying $\partial_nu(q_0) = 1$. But, since $\lambda_2$ is not a
Dirichlet eigenvalue, the space of eigenfunctions has dimension three,
whereas the multiplicity of the $\lambda_2$-Neumann eigenspace is also
three, so every eigenfunction automatically satisfies $\partial_nu(q_0) =
0$.

We note that the computation of $\beta_j$, carried out using the
recursion relation (2.11), was done using exact rational arithmetic (the
reported values are reported as decimal approximations, of course).  This
is significant because this solution of (2.11) is highly unstable.  For
example, if we take $\beta_0 = \frac 12$ and $\beta_1 = .044444444$ or
$.04444445$ (the correct value being $2/45$) and then use (2.11) for $j
\ge 2$, we find the ratio $\beta_j/\beta_{j+1}$ approaching
$-84.0799\ldots$ (this is $-5\lambda_1^D$, where $\lambda_1^D =
16.815999\ldots$ is the first Dirichlet eigenvalue).  In Section 6 we
will give an explanation for this phenomenon.

Next we will establish estimates for $\|P_{jk}\|_\infty$.  To do this we
will study the operator
$$
Af(x) = Gf(x) - (\partial_n(Gf)(q_0))P_{02}
    \tag2.28
$$
where $Gf(x) = \int G(x,y)f(y)d\mu(y)$ is the Green's operator, satisfying
$-\Delta Gf = f$ and $Gf(q_i) = 0$, $i = 0,1,2$.  Note that $A$ is a
compact linear operator, but is not self--adjoint.  Thus the spectrum of
$A$ consists of isolated eigenvalues of finite multiplicity, and zero.
Note that we have
$$
-\Delta Af = f, \ Af(q_0) = 0 \ \text{and} \ \partial_nAf(q_0) = 0.
    \tag2.29
$$
In particular, this implies
$$
AP_{jk} = -P_{(j+1)k} \ \ \text{for} \ k = 1,2.   \tag2.30
$$
Write $A_0$ for the restriction of $A$ to the $R_0$--symmetic functions,
where $R_0$ is the reflection preserving $q_0$.

\proclaim{Lemma 2.8}  (a) $f$ is an eigenfunction of $A_0$ $(A_0f =
\lambda f)$ if and only if $f$ is a symmetric
$\lambda^{-1}$--eigenfunction of $\Delta$ satisfying $f(q_0) =
\partial_nf(q_0) = 0$.  (b) $f$ is an eigenfunction of $A_0$ if and only
if $f$ is a symmetric $\lambda^{-1}$--Neumann eigenfunction of $\Delta$
satisfying $f(q_0) = 0$.  (c) The Jordan block of $A_0$ associated to any
eigenvalue is diagonal.
\endproclaim

\underbar{Proof}:  (a) By (2.29), any eigenfunction of $A$ is a
$\lambda^{-1}$--eigenfunction of $\Delta$ satisfying $f(q_0) =
\partial_nf(q_0) = 0$.  For the converse, let $v = Af-\lambda f$.  Then
$$
\Delta v = \Delta Af-\lambda\Delta f
     = \Delta(Gf-\partial_n(Gf)P_2)+f = -f+f = 0
$$
so $v$ is harmonic.  But $v$ is symmetric with $v(q_0) = \partial_nv(q_0)
= 0$, and this implies $v = 0$.

(b) The only new assertion here is that $f$ in part (a) also satisfies
$\partial_nf(q_1) = \partial_nf(q_2) = 0$.  This requires a rather
detailed knowledge of the description of eigenfunctions of $\Delta$ by
spectral decimination.  First we observe that if $|\lambda^{-1}|$ is small
enough (less than the first Dirichlet eigenvalue), then a symmetric
$\lambda^{-1}$--eigenfunction is uniquely determined by $f(q_0)$ and
$\partial_nf(q_0)$.  This implies that $f$ vanishes identically on a cell
$F^n_0(SG)$ for $n$ large enough. But an eigenfunction can vanish on a
cell only if the space of eigenfunctions has dimension greater than three,
and that happens only if $\lambda^{-1}$ is a joint Dirichlet--Neumann
eigenvalue.   That means its restriction to the graph $\Gamma_m$ for some
value of $m$ is either a 5--eigenfunction or a 6--eigenfunction.  In the
6--eigenfunction case there is nothing to prove, since all eigenfunctions
are Neumann eigenfunctions.  In the 5--eigenfunction case this is not
true, but the Neumann eigenfunctions have codimension two in the space of
all eigenfunctions.  When we impose  the $R_0$--symmetry condition the
codimension drops to one.  We know exactly what this one function looks
like (see Figure 2.2 for the case $m = 2$).  In particular, it does not
vanish identically in any small cell $F_0^m(SG)$.  Since $f$ does (and so
do all symmetric joint Dirichlet--Neumann eigenfunctions), it follows that
$f$ must be Neumann eigenfunction (in the 5--eigenfunction case it is also
a Dirichlet eigenfunction, but not necessarily in the 6--eigenfunction
case).

\bigskip
\centerline{\BoxedEPSF{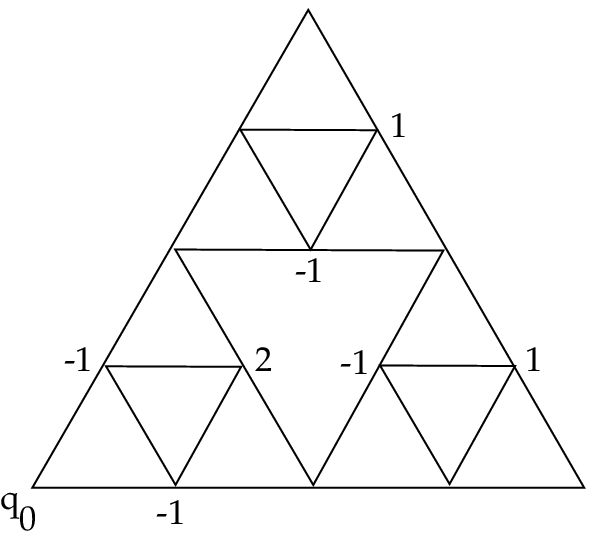 scaled 700}}

\centerline{Figure 2.2.}

\bigskip
(c) Suppose $\lambda$ is an eigenvalue of $A_0$, and $(A_0-\lambda)^2g =
0$.  Then $\lambda^{-1}$ is a Neumann eigenvalue of $\Delta$, and
$(\Delta+\lambda^{-1})^2 g = 0$.  Also $g$ is symmetric and satisfies
$g(q_0) = \partial_ng(q_0) = 0$.  By similar reasoning as before, $g$ is a
Neumann eigenfunction of $\Delta$, hence the Jordan block associated with
$\lambda$ is diagonal.  \hfill  Q.E.D.\

\bigskip
\proclaim{Theorem 2.9}  (a) For any $r < \infty$ there exists $c_r$ such
that
$$
\|P_{j1}\|_\infty \leq c_rr^{-j},   \tag2.31
$$
or more precisely
$$
\lim_{j\rightarrow\infty} \frac 1j\log \|P_{j1}\|_\infty = -\infty.
     \tag2.32
$$

(b) There exists $c$ such that
$$
\|P_{j2}\|_\infty \le c\lambda_2^{-j},   \tag2.33
$$
and
$$
\lim_{j\rightarrow\infty} (-\lambda_2)^jP_{j2} = \varphi   \tag2.34
$$
where $\varphi$ is a $\lambda_2$--Neumann eigenfunction of $\Delta$ which
is $R_0$--symmetric and vanishes on $F_0(SG)$ (a multiple of the
eigenfunction shown in Figure 2.3 on $\Gamma_1$), the limit existing
uniformly and in energy.
\endproclaim

\underbar{Proof}:  (a) Consider the norm
$$
\|f\| = (\|f\|^2_2 + \Cal E(f,f))^{1/2}
      \tag2.35
$$
and define $\Cal L_1$ and $\Cal L_2$ as the closures in this norm of the
spans of $\{P_{j1}\}$ and $\{P_{j2}\}$, respectively.  By (2.30), $A_0$
preserves both spaces.  Denote by $A_1$ and $A_2$ the restriction of $A_0$
to $\Cal L_1$ and $\Cal L_2$. We claim $\sigma(A_1) = \{0\}$.  Indeed,
otherwise $A_1$ would have to have a nonzero eigenvalue $\lambda$ because
$A_1$ is compact.  Since this would also be an eigenvalue of $A_0$, by
Lemma 2.8 $\lambda^{-1}$ would have to be a Neumann eigenvalue of
$\Delta$.  So $\lambda > 0$, and we may choose it to be the largest
eigenvalue of
$A_1$.  Then $\lambda^{-j}A^j_1$ converges to a projection (not
necessarily orthogonal) $B_\lambda$ onto the finite dimensional
$\lambda$--eigenspace of $A_1$.  Note that $B_\lambda P_{01}$ cannot be
the zero function, because that would imply $B_\lambda P_{j1} = 0$ for
all $j$, contradicting the fact that $B_\lambda$ is nonzero.  But then
$\lambda^{-j}A^j_1P_{01} = \lambda^{-j}P_{j1}$ would converge to a nonzero
eigenfunction of $A_1$.  By Theorem 2.7 this eigenfunction would vanish at
$q_1$ and $q_2$, and of course it vanishes at $q_0$, since $P_{j1}$ does
for $j \ge 1$.  So it would have to be a joint Dirichlet--Neumann
eigenfunction of $\Delta$.  But Theorem 4.3.6 of [Ki2] asserts that all
$P_{jk}$ are orthogonal to all joint Dirichlet--Neumann eigenfunctions.

Thus we have shown that $\sigma(A_1) = \{0\}$, so the spectral radius of
$A_1$ is zero,
$$
\lim_{j\rightarrow\infty} \|A^j_1\|^{1/j} = 0.
$$
Applying this to $P_{01}$ we obtain (2.32) (the norm (2.35) dominates the
$L^\infty$ norm), which implies (2.31).

(b) The result of Kigami used above moreover says that $\Cal L = \Cal L_1
\oplus \Cal L_2$ contains all $R_0$--symmetric Neumann eigenfunctions of
$\Delta$ that are othogonal to all joint Dirichlet--Neumann eigenfunctions
(note that Kigami uses the $L^2$ norm rather than (2.35), but the same
argument applies).  In particular, it contains the
$\lambda_2$--eigenfunction shown in Figure 2.3 (this is a Neumann
eigenfunction, so it is orthogonal to all Neumann eigenfunctions with
different eigenvalues, and there are \linebreak
no joint Dirichlet--Neumann
eigenfunctions with the same eigenvalue).  By Lemma 2.8 and the explicit
description of Neumann eigenfunctions, $\lambda_2^{-1}$ is the largest
eigenvalue of $A_0$, and $\varphi$ spans this multiplicity one eigenspace.
Thus, as before, $\lambda^j_2A^j$ converges to a one--dimensional
projection operator $B_{\lambda_2^{-1}}$, and $B_{\lambda_2^{-1}}P_{01} =
0$.  That means $B_{\lambda_2^{-1}} P_{02} \neq 0$, for otherwise
$B_{\lambda_2^{-1}} = 0$.  So
$$
\lim_{j\rightarrow\infty} (-\lambda_2)^jP_{j2}
       = \lim_{j\rightarrow\infty} \lambda^j_2A^jP_{02}
     = B_{\lambda_2^{-1}} P_{02}
$$
which is (2.34).  This implies (2.33).  \hfill  Q.E.D.

\bigskip
The estimate (2.33) is sharp, but (2.32) falls short of what we would have
if we knew $\|P_{j1}\|_\infty = \alpha_j$, in view of (2.27).  One
approach to establish this would be to prove the following conjecture:

\

\centerline{\BoxedEPSF{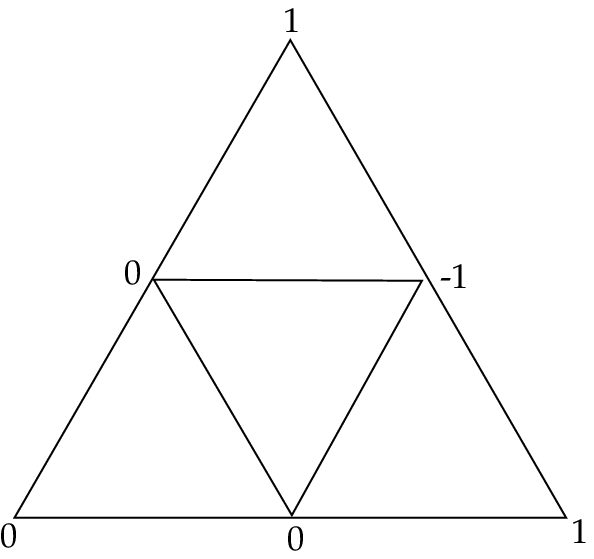 scaled 700}}

\medskip
\centerline{Figure 2.3}.

\bigskip


\proclaim{Conjecture 2.10}  For all $x \neq q_0$ and all $j$,
$$
P_{j1}(x) > 0.    \tag2.36
$$
\endproclaim

We have numerical evidence for this conjecture for moderate values of
$j$.  To show that (2.36) implies $\|P_{j1}\|_\infty = \alpha_j$ is easy
using the following well-known fact (we provide a proof since it does not
appear explicitly in the literature).

\proclaim{Proposition 2.11}  If $u \in$ dom $\Delta$, $\Delta u(x_0) > 0$
and $x_0$ is not a boundary point, then $u$ does not achieve its maximum
value at $x_0$.
\endproclaim

\underbar{Proof}:  If $x_0$ is a vertex in $V_*$ the result follows
immediately from the pointwise definition of $\Delta u(x_0)$.  If not,
then we can find a cell $F_wK$ such that $x_0$ is in the interior of
$F_wK$ and $\Delta u > 0$ on $F_wK$.  Let $v = u \circ F_w$.  Then
$\Delta v > 0$, and we have
$$
v(x) = h(x) - \int_K G(x,y)\Delta v(y)dy
$$
where $G$ is the Dirichlet Green's function and $h(x)$ is the harmonic
function with the same boundary values as $v(x)$.  Since the Green's
function is positive in the interior, we have $v(x) < h(x)$ in the
interior.  Since $h$ attains its maximum on the boundary, it follows
that $v$ cannot attain its maximum in the interior, so $u(x_0)$ is not a
maximum.    \hfill  Q.E.D.

\bigskip
Next we study the normal and tangential derivatives  of monomials at
boundary points.

\proclaim{Theorem 2.12}  We have initial values $n_0 = 0$, $t_0 = -1/2$,
and recursion relations
$$
n_j = \frac{5^j+1}{2} \alpha_j + 2 \sum^{j-1}_{\ell=0} n_\ell
        \beta_{j-\ell} \quad \text{for} \quad j \ge 1,
      \tag2.37
$$
$$
t_j = \beta_j - 6 \sum^{j-1}_{\ell=0} \alpha_{j+1-\ell} t_\ell \quad
        \text{for} \quad j \ge 1.
     \tag2.38
$$
Moreover, we have
$$
\partial_n P_{j2}(q_1) = \partial_n P_{j2}(q_2)
   = \cases  \frac 12 - \alpha_0  &\text{if} \quad j = 0 \\
           -\alpha_j &\text{if} \quad j \ge 1
     \endcases
    \tag2.39
$$
$$
\partial_n P_{j3}(q_1) = -\partial_n P_{j3}(q_2) = 3n_{j+1}
    \tag2.40
$$
$$
\partial_TP_{j1}(q_1) = -\partial_TP_{j1}(q_2)
    = \cases \frac 16  &\text{if} \quad j = 1 \\
               0  &\text{if} \quad j \neq 1
      \endcases
     \tag2.41
$$
$$
\partial_TP_{j3}(q_1) = -\partial_TP_{j3}(q_2)
    = \cases -\frac 12  &\text{if} \quad j = 0 \\
               0  &\text{if} \quad j \ge 1
      \endcases.
     \tag2.42
$$
\endproclaim

\underbar{Proof}:  As in the proof of Theorem 2.3 we introduce matrices
$n$, $\wt n$ and $t$, where $\wt n_j = \partial_n P_{j2}(q_1)$.   When we
evaluate the normal derivatives on both sides of (2.14) at
$q_1$, we see that
$$    n_j = b_{j-1} +\sum_{l=0}^j \alpha_{j-l}(a_{l-1}+b_{l-1})\qquad
     \text{for all} \ j,
$$
or in matrix notations
$$
     n = B + \alpha (A+B).
$$
Using (2.15) this yields
$$
n = \frac 12 \beta^{-1} (I+2\alpha)(I-\alpha)
    = \frac 14 \beta^{-1}(2I-\tau(\alpha)-\alpha)
   \tag2.43
$$
which implies (2.37).

By the same reasoning
$$
\wt n_j = \sum^j_{l=0} \beta_{j-l}(a_{l-1}+b_{l-1}) \ \text{for all} \
        j.
$$
Then
$$
\wt n = \beta(A+B)
$$
and hence by (2.15) we obtain
$$
\wt n = \frac 12 I - \alpha,   \tag2.44
$$
which implies (2.39).

Finally, the same reasoning shows
$$
t_j = \sum^j_{l=0} \beta_{j-\ell}T_l \ \text{for all} \ j,
$$
where $T_l = \partial_Tf_{l2}(q_1)$.  Now $P_{j3} =
\sum^j_{l=0}\gamma_{j-l}(f_{l1}-f_{l2})$, so taking tangential derivatives
at $q_0$ we get
$$
2 \sum^j_{l=0} \gamma_{j-l}T_l = \partial_TP_{j3}(q_0)
    = \cases 1,  &\text{if} \ j = 0; \\
             0, &\text{otherwise}.
      \endcases
$$
In matrix notations these become
$$
t = \beta T
$$
$$
\gamma T = \frac 12 I.
$$
Together we have
$$
\beta = 2\gamma t = 6\sigma(\alpha)t,  \tag2.45
$$
where the last equality follows form (2.12).

This proves (2.38).  The initial values of $n_0$, $\wt n_0$ and $t_0$ are
easy to check.

Note that the skew-symmetry implies $\partial_TP_{j3}(q_1) =
\partial_TP_{j3}(q_2)$, so (2.2) implies $\partial_TP_{j3}(q_0) +
2\partial_TP_{j3}(q_1) = 0$, which yields (2.42).  Then (2.41) follows
from (2.19) and (2.42), and similarly (2.19) implies (2.40).  \hfill
Q.E.D.

\bigskip

\proclaim{Theorem 2.13}  For any $r < \infty$ there exists $c_r$ such
that, for all $j \ge 1$,
$$
|n_j|\le c_rr^{-j}.   \tag2.46
$$
Also
$$
|t_j| \le c\lambda_2^{-j}.   \tag2.47
$$
\endproclaim

\underbar{Proof}:  From the Gauss--Green formula we have
$$
\int \Delta ud\mu = \sum^2_{i= 0} \partial_nu(q_i).
$$
We apply this to $u = P^{(0)}_{j1}$, noting that
$\partial_nP^{(0)}_{j1}(q_0) = 0$ and $\partial_nP^{(0)}_{j1}(q_1) =
\partial_nP^{(0)}_{j1}(q_2) = n_j$.  It follows that
$$
n_j = \frac 12 \int P^{(0)}_{(j-1)1} d\mu,
    \tag2.48
$$
and (2.46) follows from (2.31).

Similarly, (2.47) will follow from (2.33) and the estimate
$$
|\partial_Tu(q_i)|
      \le c(\|u\|_\infty + \|\Delta u\|_\infty + \|\Delta^2u\|_\infty).
    \tag2.49
$$
In [S3] it is shown that $\partial_Tu(q_i)$ exists if $u \in \dom\Delta$
and $\Delta u$ satisfies a H\"older condition, and (2.49) is just a
quantitative version of this fact.  For the convenience of the reader we
outline the argument. For simplicity take $i = 0$.  Let $g_m$ (see Figure
2.4 for $m = 2$) denote the level $m$ piecewise harmonic function
satisfying $g_m(q_0) = 0$ and $g_m(F^k_0q_1) = 3^k$ and $g_m(F^k_0q_2) =
-3^k$ for all $k \le m$.  Then

\bigskip

\centerline{\BoxedEPSF{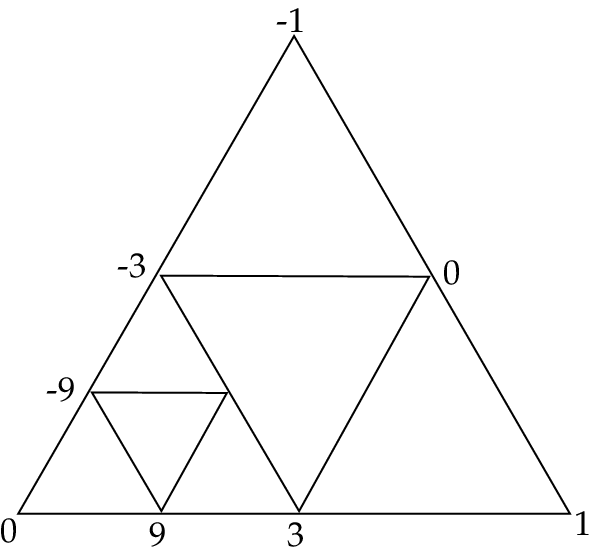 scaled 700}}


\bigskip
\centerline{Figure 2.4}.

\medskip
$$
\int g_m\Delta u d\mu
     = \frac{14}{3} 5^m(u(F^m_0q_1) - u(F^m_0q_2))
       - 5(u(q_1)-u(q_2))
    \tag2.50
$$
by the Gauss--Green formula, since the sum of the normal derivatives of
$g_m$ at $F^m_0q_1$ is $(14/3)5^m$ (there are no terms involving normal
derivatives of $u$ at $F^m_0q_i$ because $u$ satisfies matching
conditions).  Let $u_1 = \Delta u$.  Note that $g_m$ is odd, so only the
odd part of $u_1$ contributes to the integral in (2.50).  So (2.49) will
follow from (2.50) and the estimate
$$
\Big|\int g_m(u_1-u_1\circ R_0)d\mu\Big|
   \le c(\|u_1\|_\infty + \|\Delta u_1\|_\infty).
     \tag2.52
$$
But (2.52) is routine, because on the cells $F^k_0F_1(SG)$ and
$F^k_0F_2(SG)$ $(0 \le k \le m)$ of measure $3^{-k-1}$, the function $g_m$
is of size $3^k$, and $u_1-u_1\circ R_0$ can be estimated by  $(\frac
35)^k \|\Delta u_1\|_\infty$.   \hfill  Q.E.D.

\bigskip
In Table 5.2 we display the results of solving the recursion relations for
$n_j$ and $t_j$.  The data suggests that $(-\lambda_2)^jt_j$ converges, in
fact quite a bit faster than for $\beta_j$, and
$\ds\lim_{j\rightarrow\infty} \beta_j/t_{j+1} = 9$.  Moreover $n_j$ is
always positive and satisfies
$$
n_j \le cj\alpha_j.   \tag2.53
$$
If Conjecture 2.10 holds, then $\|P_{(j-1)1}\|_\infty = \alpha_{j-1}$ so
(2.48) implies $n_j \le \frac 12 \alpha_{j-1}$, which is only slightly
weaker than (2.53).


\bigskip
$$
\matrix
j   &n_j  &t_j  &\frac{n_j}{j\alpha_j}   &(-\lambda_2)^jt_j
                & \frac{\beta_j}{t_{j+1}} \\ \ \\
0    &  0   &-.50000000      &\infty  &-.50000000    &18. \\
1   &.50000000    &-.027777778    &3.    &3.7658925
             &-432. \\
2   &.027777778   &.00010288066   &2.5000000   &1.8909261
             &1439.0526 \\
3  &.00041152263   &-.70062097\  10^{-6}   &2.2222222   &1.7457996
           &-1679.0103 \\
4   &.27287343 \ 10^{-5}   &.50952342\  10^{-8}   &2.0555556
            &1.7212575    &1027.9833 \\
5   &.98752993 \ 10^{-8}   &-.37481616\  10^{-10}   &1.9341564
            &1.7166051   &-356.40392 \\
6   &.22167060 \ 10^{-10}   &.27632364 \ 10^{-12}   &1.8405958
           &1.7156968    &94.889369 \\
7   &.33533009 \ 10^{-13}   &-.20379909 \ 10^{-14}   &1.7656562
           &1.7155176     &-5.1358463 \\
8   &.36203261 \ 10^{-16}   &.15032210\  10^{-16}   &1.7035627
           &1.7154821   &10.708638 \\
9   &.29106143 \ 10^{-19}   &-.11087934 \ 10^{-18}   &1.6507112
           &1.7154750   &8.8428158 \\
10   &.18012308\  10^{-22}   &.81786167 10^{-21}   &1.6048457
              &1.7154736   &9.0113344 \\
11   &.88115370\  10^{-26}   &-.60326673\  10^{-23}   &1.5644762
             &1.7154734    &8.9993459  \\
12   &.34823920\  10^{-29}   &.44497842 \ 10^{-25}   &1.5285491
                &1.7154734   &9.0000311 \\
13   &.11321107\  10^{-32}   &-.32822264 \ 10^{-27}  &1.4962768
          &1.7154734     &8.9999988 \\
14   &.30738762\  10^{-36}   &.24210186 \ 10^{-29}   &1.4670507
           &1.7154734    &9.0000000 \\
15   &.70615767 \ 10^{-40}   &-.17857790 \ 10^{-31}   &1.4403911
          &1.7154735   &9.0000000 \\
16   &.13880322\  10^{-43}   &.13172169 \ 10^{-33}  &1.4159159
           &1.7154735    &9.0000000  \\
17  &.23573795 \ 10^{-47}   &-.97159864 \ 10^{-36}  &1.3933188
            &1.7154736    &9.0000000  \\
18   &.34893132\  10^{-51}   &.71666548\  10^{-38}   &1.3723521
          &1.7154736    &9.0000000 \\
19   &.45359082 \ 10^{-55}   &-.52862303\  10^{-40}   &1.3528138
          &1.7154736    &9.0000000  \\
20   &.52141937 \ 10^{-59}   &.38992014\  10^{-42}  &1.3345373
            &1.7154737   &9.0000000
\endmatrix
$$

\centerline{Table 2.2.}
\bigskip

\medskip
\noindent
  We also have found that the recursion relation for $n_j$ is
unstable, and any slight perturbation produces a decay rate
$O((\lambda_1^D)^{-j})$, which is even slower than the decay rate for
$\beta_j$ and $t_j$. Also a slight perturbation of the $t_j$ recursion
relation produces a decay rate of $O((\lambda^D_2)^{-j})$. We will
explain this in Section 6.

Next we describe the change of basis formula to pass between
$\{P_{jk}^{(\ell)}\}$ for different values of $\ell$, an immediate
consequence of Theorem 2.12.

\proclaim{Corollary 2.14}  We have
$$
\pmatrix  P^{(\ell)}_{j1} \\ P^{(\ell)}_{j2} \\ P^{(\ell)}_{j3}
\endpmatrix
= \sum^j_{k=0} M_{j-k}
\pmatrix  P^{(\ell+1)}_{k1} \\ P^{(\ell+1)}_{k2} \\
         P^{(\ell+1)}_{k3}
\endpmatrix
    \tag2.54
$$
for matrices $M_j$ given by
$$
\cases  M_j =
    \pmatrix \alpha_j  &n_j  &0 \\ \beta_j  &-\alpha_j  &t_j \\
           3\alpha_{j+1}  &3n_{j+1}  &0
     \endpmatrix \quad \text{for} \quad j \ge 2 \\
M_1 = \pmatrix  \alpha_1  &n_1  &\frac 16 \\  \beta_1  &-\alpha_1 &t_1 \\
            3\alpha_2  &3n_2  &0
      \endpmatrix  \quad
M_0 =  \pmatrix \alpha_0  &n_0  &0 \\ \beta_0  &\frac 12-\alpha_0 &t_0\\
           3\alpha_1  &3n_1  &-\frac 12
    \endpmatrix.
\endcases   \tag2.55
$$
Similarly
$$
\pmatrix  P^{(\ell)}_{j1} \\ P^{(\ell)}_{j2} \\ P^{(\ell)}_{j3}
\endpmatrix
= \sum^j_{k=0} \wt M_{j-k}
     \pmatrix  P^{(\ell-1)}_{k1} \\ P^{(\ell-1)}_{k2} \\
          P^{(\ell-1)}_{k3}
     \endpmatrix
   \tag2.56
$$
for
$$
\cases
\wt M_j =
     \pmatrix \alpha_j  &n_j  &0 \\ \beta_j  &-\alpha_j  &-t_j \\
                -3\alpha_{j+1}  &-3n_{j+1}  &0
      \endpmatrix \quad \text{for} \quad j \ge 2 \\
\wt M_1 =
   \pmatrix  \alpha_1  &n_1  &-\frac 16 \\ \beta_1  &-\alpha_1  &-t_1 \\
            -3\alpha_2  &-3n_2  &0
    \endpmatrix\quad
\wt M_0 =
    \pmatrix  \alpha_0  &n_0  &0 \\  \beta_0  &\frac 12-\alpha_0  &-t_0 \\
           -3\alpha_1  &-3n_1  &-\frac 12
     \endpmatrix.
\endcases
   \tag2.57
$$
\endproclaim

\vskip.4in
\subheading{\S3.  Power series}

\medskip
A formal power series about $q_\ell$ is an expression of the form
$$
\sum^3_{k=1} \ \sum^\infty_{j=0} c_{jk} P^{(\ell)}_{jk}(x).
     \tag3.1
$$
We call $\{c_{jk}\}$ the {\it coefficients}, and we seek growth
conditions on the coefficients that will make (3.1) converge nicely.

\proclaim{Theorem 3.1}  If the coefficients satisfy
$$
|c_{j1}| \ \text{and} \ |c_{j3}| = O((j!)^r)
         \ \text{for some} \ r < \log 5/\log 2,
    \tag3.2
$$
and
$$
|c_{j2}| = O(R^j) \ \text{for some} \ R < \lambda_2
    \tag3.3
$$
then (3.1) converges uniformly and absolutely to a function $u \in$
dom\,$(\Delta^\infty)$, and (3.1) may be ``differentiated term-by-term'',
$$
\Delta^nu(x) = \sum^3_{k=1} \ \sum^\infty_{j=n} c_{jk}
       P^{(\ell)}_{(j-n)k}(x).
   \tag3.4
$$
Moreover, the coefficients are given by the infinite jet of $u$ at
$q_\ell$:
$$
\cases
   c_{j1} = \Delta^j u(q_\ell) \\
  c_{j2} = \partial_n\Delta^ju(q_\ell) \\
  c_{j3} = \partial_T \Delta^ju(q_\ell).
\endcases
   \tag3.5
$$
\endproclaim

\underbar{Proof}:  The estimates in Theorem 2.9 conspire with the
growth rates (3.2) and (3.3) to make (3.1) converge uniformly and
absolutely.  Call the limit $u$.  Note that the right side (3.4) is also
a formal power series, in fact
$$
\sum^3_{k=1} \ \sum^\infty_{j=0} c_{(j+n)k} P^{(\ell)}_{jk}(x)
$$
whose coefficients also satisfy the growth rate conditions (3.2) and
(3.4).  So the right side of (3.4) converges uniformly and absolutely.
By terminating the sums at $j = N$ and letting $N \rightarrow \infty$ we
obtain the equality in (3.4) by a routine argument using the Green's
function [Ki2].

It suffices to prove the jet formulas (3.5) when $j = 0$ in view of
(3.4), and for this it suffices to show that if $c_{01} = c_{02} = c_{03}
= 0$ then $u(q_\ell) = \partial_nu(q_\ell) = \partial_Tu(q_\ell) = 0$.
Of course $u(q_\ell) = 0$ directly from (3.1).  For simplicity put $\ell
= 0$.  Then (since $u(q_0) = 0$)
$$
\partial_nu(q_0) = -\lim_{m\rightarrow\infty}\Big(\frac 53\Big)^m
         (u(F_0^mq_1) + u(F_0^mq_2)).
$$
But we have
$$
u(F^m_0x) = \sum^\infty_{j=1} c_{j1}5^{-mj} P_{j1}(x)
       + c_{j2}\Big(\frac 35 5^{-j}\Big)^m P_{j2}(x)
       + c_{j3} 5^{-m(j+1)} P_{j3}(x).
    \tag3.6
$$
Using the estimates for the coefficients and monomials we see that
$$
u(F^m_0x) = O(5^{-m}),   \tag3.7
$$
and this suffices to prove $\partial_nu(q_0) = 0$.  This by itself does
not suffice for the tangential derivative, which has a factor of $5^m$.
However, for the tangential derivative we can restrict attention to the
skew-symmetric part
$$
\wt u(x) = \frac 12 (u(x)-u(\rho_0x))
      = \sum^\infty_{j=1} c_{j3}P_{j3}(x),
    \tag3.8
$$
so the analog of (3.6) shows
$$
\wt u(F^m_0x) = O(5^{-2m}),   \tag3.9
$$
which implies $\partial_T u(q_0) = 0$.  \hfill  Q.E.D.

\bigskip
As a corollary of the proof we can characterize rates of vanishing of
power series.

\medskip
\noindent
\underbar{Definition 3.2}:  A function $f$ is said to {\it vanish to
order} $r$ (any positive real) at $q_\ell$ provided
$$
\|f\circ F^m_\ell\|_\infty = O(5^{-mr}).   \tag3.10
$$
If (3.10) holds for all $r$ then we say $f$ {\it vanishes to infinite
order} at $q_\ell$.

\proclaim{Corollary 3.3}  If $u$ is represented
by a power series (3.1) with coefficients satisfying growth conditions
(3.2) and (3.3), then $u$ vanishes to order $N$ (a positive integer) at
$q_\ell$ if and only if $c_{jk} = 0$ for all $j < N$.  In that case
$\Delta^\ell u$ vanishes to order $N-\ell$ for all $\ell < N$.  Moreover,
the odd part $\wt u$ vanishes to order $N+1$.  In particular, if $u$ is
not identically zero then it cannot vanish to infinite order.
\endproclaim

Next we consider rearrangement of power series, moving from one boundary
point $q_\ell$ to another.  It turns out that we need to make stronger
assumptions on the coefficients, requiring $c_{j1}$ and $c_{j3}$ to
satisfy the same exponential growth rate as $c_{j2}$.

\proclaim{Theorem 3.4}   Suppose the
coefficients of a power series (3.1) about one boundary point $q_\ell$
satisfy
$$
|c_{jk}| = O(R^j) \ \text{for some} \ R < \lambda_2, \
          k = 1,2,3.
      \tag3.11
$$
Then the function may also be represented by power series about the other
boundary points with coefficients also satisfying (3.11).  More
precisely, the coefficients at $q_{\ell+1}$ are given by
$$
(c'_{j'1}\ c'_{j'2}\ c'_{j'3}) = \sum^\infty_{j=0}
        (c_{(j+j')1} \ c_{(j+j')2}\  c_{(j+j')3})M_j
      \tag3.12
$$
and similarly at $q_{\ell-1}$ with $M_j$ replaced by $\wt M_j$ (see
(2.55) and (2.57)).
\endproclaim

\underbar{Proof}:  The key observation is that the right side of (3.12)
converges absolutely and the new coefficients again satisfy (3.11) (in
fact with the same value of $R$) because the entries in $M_j$ are
$O(\lambda_2^{-j})$ by Theorem 2.13.  Of course (3.11) is exactly what
we get if we substitute (2.54) into (3.1) and interchange the order of
summation, which is easily justified using the estimates of Theorem
2.9.  \hfill  Q.E.D.

\bigskip
Note that we could not allow slower growth rates like (3.2) for the
$c_{j1}$ and $c_{j3}$ coefficients and still rearrange, because the
second column of $M_j$ has positive entries.  In \S6 we will present an
example to show that rearrangement fails when $c_{j1} = O(\lambda^j_2)$.
However, condition (3.11) is not sharp.  We could replace it by
$$
\sum^\infty_{j=0} \lambda_2^{-j}|c_{jk}| < \infty,
    \tag3.13
$$
and the rearranged coefficients would satisfy the same growth condition.
However, not all subsequent results would be valid under this hypothesis.

\medskip
\noindent
\underbar{Definition 3.5}:  An {\it entire analytic function} is a
function given by a power series (3.1) with coefficients satisfying
(3.11).

We can also consider local power series expansions on any cell $F_w(SG)$
with respect to a boundary point $F_wq_\ell$ of the cell, namely
$$
\split
\sum^\infty_{j=0} \Big(5^{-mj}c_{j1}&P_{j1}^{(\ell)}(F_w^{-1}x)
          + \Big(\frac 35 5^{-j}\Big)^{m} c_{j2}P_{j2}^{(\ell)}
             (F_w^{-1}x) \\
     &+ 5^{-(j+1)m} c_{j3} P_{j3}^{(\ell)} (F_w^{-1}x)\Big)
\endsplit   \tag3.14
$$
where $m = |w|$.

\proclaim{Theorem 3.6}  An entire
analytic function has a local power series expansion (3.14) for any $w$
and $\ell$ with coefficients satisfying (3.11).  Conversely, suppose
$u(x)$ is a function defined on $F_w(SG)$ given by a local power series
expansion (3.14) with coefficients satisfying (3.11).  Then $u$ has a
unique extension to an entire analytic function.
\endproclaim

\underbar{Proof}:  Suppose first that $m = 1$, say $w = (0)$.  If $\ell =
0$ then the local and global power series are identical, with identical
coefficients.  Moreover, $u \circ F_w$ is an entire analytic function
with coefficients satisfying (3.11) (in fact with $R < \lambda_2/5$).
The rearrangement for $u \circ F_w$ about $q_1$ and $q_2$ guaranteed by
Theorem 3.4 gives the local power series of $u$ in $F_0(SG)$ about
$F_0q_1$ and $F_0q_2$, with the same coefficient estimates.  We may then
iterate this argument to get local power series about any boundary point
in any cell.

Conversely, suppose $u$ is given in $F_w(SG)$ by a local power series
about $F_wq_\ell$, with coefficients satisfying (3.11).  Write $w =
(w',w_m)$ with $|w'| = m-1$.  If $w_m \neq \ell$ then use Theorem 3.4 to
rearrange the power series of $u\circ F_w$ about $q_{w_m}$.  So we end up
with a local power series of $u$ about $F_{w'}F_\ell q_\ell$ in the cell
$F_{w'}F_\ell(SG)$.  But $F_{w'}F_\ell q_\ell = F_{w'}q_\ell$ and the
power series makes sense in the cell $F_{w'}(SG)$.  Use this power series
to extend the definition of $u$.  By iterating the argument, we obtain
the desired extension.  Note that the estimates (3.11) on the
coefficients are reproduced in each extension or rearrangement step.  It
is clear that the extension is unique because the rearranged coefficients
are determined by (3.12).  \hfill Q.E.D.

\bigskip
By the same reasoning, if a local power series has coefficients satisfying
$$
c_{jk} = O(R^j) \ \text{for some} \ R < 5^{m_0}\lambda_2,
    \tag3.15
$$
then the function can be also represented by a power series on a level
$m_0$ cell.  One might hope that this ``analytic continuation'' might
extend somewhat beyond the cell, with the domain of analyticity growing
as $R$ decreases toward $5^{m_0-1}\lambda_2$.  However, the experimental
evidence we have seen does not support this at all.  On the contrary, we
will see in \S6 that there are power series (3.1) with coefficients
$O(\lambda^j_2)$ where we have divergence outside $F_\ell(SG)$.  We might
describe this as a ``quantized radius of convergence.''  Of course, this
does not rule out a different type of behavior for special classes of
power series.

\proclaim{Theorem 3.7}  An entire
analytic function satisfies the estimate
$$
\|\Delta^nu\|_\infty = O(R^n) \ \text{for some} \ R < \lambda_2.
     \tag3.16
$$
\endproclaim

\underbar{Proof}:  We have
$$
\Delta^nu = \sum^3_{k=1} \ \sum^\infty_{j=n} c_{jk} P^{(\ell)}_{(j-n)k}
$$
so
$$
\|\Delta^nu\|_\infty
     \le M \sum^3_{k=1} \ \sum^\infty_{j=n}
              R^j\|P^{(\ell)}_{(j-n)k}\|_\infty
    \le M \sum^\infty_{j=n} R^j\lambda_2^{n-j} = O(R^n)
$$
for $R$ in (3.11).    \hfill Q.E.D.

\bigskip
The condition (3.16) obviously implies the same estimate in $L^2$ norm:
$$
\|\Delta^nu\|_2 = O(R^n) \ \text{for some} \ R < \lambda_2.
    \tag3.17
$$
But conversely, (3.17) implies (3.16), because
$\|f\|_\infty \le c(\|f\|_2 + \|\Delta f\|_2)$.  The estimate (3.17) is
technically more convenient, since we can compute $L^2$ norms exactly
from eigenfunction expansions.

It follows immediately from the definition that an eigenfunction of
$\Delta$ is an entire analytic function if and only if the eigenvalue
satisfies $|\lambda| < \lambda_2$.  Theorem 3.7 shows us that many other
functions that we might believe to be entire analytic functions are not.
Indeed, suppose $u$ is represented by a Dirichlet (or Neumann)
eigenfunction expansion
$$
u(x) = \sum^\infty_{k=1} a_k\varphi_k(x)
     \tag3.18
$$
where $\{\varphi_k\}$ is an orthonormal basis of Dirichlet (or Neumann)
eigenfunctions.  If the coefficients are rapidly decreasing,
$$
a_k = O(k^{-n}) \ \text{for all} \ n,   \tag3.19
$$
then we may differentiate term-by-term,
$$
\Delta^nu(x) = \sum^\infty_{k=1} (\lambda^D_k)^n a_k\varphi_k(x).
     \tag3.20
$$
It follows that
$$
\|\Delta^nu\|_2 = \Big(\sum^\infty_{k=1} (\lambda^D_k)^{2n} |a_k|^2
        \Big)^{1/2}.
      \tag3.21
$$
If (3.18) is non-trivial in the sense that an infinite number of
coefficients are non-zero, then not only does (3.17) fail to hold, but
the estimate cannot hold for any finite $R$.  So $u$ cannot be represented
by a local power series with (3.14) holding on any cell.  In particular
this applies to the heat kernel.

This observation stands in striking contrast to the situation on the unit
interval, where analyticity properties of a function may be characterized
by decay properties of the coefficients of its Fourier series expansion.

\bigskip
\subheading{\S4.  Characterization of analytic functions}

\medskip
The main purpose of this section is to prove the following theorem.

\proclaim{Theorem 4.1}  $u$ is an
entire analytic function if and only if $u \in dom (\Delta^\infty)$
and (3.16) (or equivalently (3.17)) holds.
\endproclaim

We first consider the case when $u$ is even with respect to $\rho_0$.  In
that case we would like a Taylor expansion with remainder about $q_0$,
$$
u(x) = T_ku(x) + R_k(x)
     \tag4.1
$$
for
$$
T_ku(x) = \sum^{k-1}_{j=0} \Delta^ju(q_0) P_{j1}(x)
        + (\partial_n\Delta^ju(q_0))P_{j2}(x)
    \tag4.2
$$
and $R_k(x)$ the remainder term.  While we can use (4.1) to define the
remainder, to be useful we need some explicit expression for it.  We are
only able to do this for $x = q_1$ (or $q_2$).

\proclaim{Lemma 4.2}  Let $v_k$ be a function in $\Cal H_k$ that is even
with respect to $\rho_0$ satisfying
$$
\Delta^jv_k(q_1) = 0 \ \text{for} \ j \le k-1   \tag4.3
$$
$$
\partial_n\Delta^jv_k(q_1)
   = \cases  0  &\text{for} \ j \le k-2 \\
            -\frac 12  &\text{for} \ j = k-1.
     \endcases
    \tag4.4
$$
Then
$$
R_k(q_1) = R_k(q_2) = \int_{SG} v_k\Delta^kud\mu
     \tag4.5
$$
for even functions $u \in dom(\Delta^k)$.
\endproclaim

\underbar{Proof}:  Note that $\Delta^ku = \Delta^k(u-T_ku) =
\Delta^kR_k$.  We apply the Gauss-Green formula $k$ times to obtain
$$
\split
\int v_k\Delta^ku d \mu
   &= \int v_k\Delta^kR_k d\mu \\
    &= 2 \sum^{k-1}_{j=0} \big(\Delta^jv_k(q_1) \partial_n\Delta^{k-j-1}
           R_k(q_1)
      - \partial_n\Delta^jv_k(q_1)\Delta^{k-j-1} R_k(q_1)\big)
\endsplit
$$
since $\Delta^{k-j-1}R_k(q_0) = \partial_n\Delta^{k-j-1}(q_0) = 0$.  By
(4.3) and (4.4) all terms vanish except when $j = k-1$ and we obtain
exactly $R_k(q_1)$.   \hfill  Q.E.D.

\bigskip
\proclaim{Lemma 4.3}  The function
$$
v_k = \sum^{k-1}_{\ell=0} (-\beta_{k-\ell-1} P_{\ell 1}^{(0)}
            + \alpha_{k-\ell-1}P_{\ell 2}^{(0)})
    \tag4.6
$$
satisfies the conditions of Lemma 4.2.
\endproclaim

\underbar{Proof}:  Clearly $v_k \in \Cal H_k$ and is even.  Since
$$
\Delta^jv_k = \sum^{k-j-1}_{\ell=0} -\beta_{k-j-1-\ell} P_{\ell 1}
         + \alpha_{k-j-1-\ell} P_{\ell 2}
$$
we obtain
$$
\Delta^jv_k(q_1) = \sum^{k-j-1}_{\ell =0} (-\beta_{k-j-1-\ell}\alpha_\ell
          + \alpha_{k-j-1-\ell} \beta_\ell) = 0
$$
which is (4.3).  Similarly
$$
\partial_n\Delta^jv_k(q_1)
     = \sum^{k-j-1}_{\ell =0} \Big(-\beta_{k-j-1-\ell} n_\ell
        - \sum^{k-j-1}_{\ell =0} \alpha_{k-j-1-\ell} \alpha_\ell\Big)
     + \frac 12 \alpha_{k-j-1}
$$
by (2.35). When $j = k-1$ this is just
$$
\partial_n\Delta^{k-1}v_k(q_1)
      = \beta_0n_0 - \alpha^2_0 + \frac 12 \alpha_0 = -\frac 12.
$$
For $j \le k-2$ we have
$$
\sum^{k-j-1}_{\ell =0} \beta_{k-j-1-\ell} n_\ell
    = - \Big(\frac{5^{k-j-1}+1}{4}\big)\alpha_{k-j-1}
$$
by (2.33), and
$$
\sum^{k-j-1}_{\ell =0} \alpha_{k-j-1-\ell} \alpha_\ell
      = \Big(\frac{5^{k-j-1}+3}{4}\Big)\alpha_{k-j-1}
$$
by (2.9) (this uses $k-j-1 \ge 1$).  Thus $\partial_n\Delta^jv_k(q_1) =
0$, proving (4.4).   \hfill  Q.E.D.

\bigskip
\proclaim{Lemma 4.4}   If $u$ is an even
function in $dom(\Delta^k)$ satisfying (3.16), and $\wt u$ is the entire
analytic function whose expansion about $q_0$ has coefficients $c_{j1} =
\Delta^ju(q_0)$, $c_{j2} = \partial_n \Delta^ju(q_0)$ and $c_{j3} = 0$,
then $u(q_1) = \wt u(q_1)$ and $u(q_2) = \wt u(q_2)$.
\endproclaim

\underbar{Proof}:  First we observe that (3.16) implies the coefficients
of $\wt u$ satisfy (3.11). This is obvious for $c_{k1}$ and $c_{k3}$, but
it follows for $c_{k2}$ because $\partial_nf(q_0) = \int hfd\mu$ for a
fixed harmonic function $h$.  Now apply Lemma 4.2 to the function $u-\wt
u$ to obtain
$$
|u(q_1) - \wt u(q_1)|
      = \Big|\int v_k\Delta^k(u-\wt u)d\mu\Big|
     \le cR^k\|v_k\|_\infty.
$$
But we easily obtain $\|v_k\|_\infty = O(\lambda^{-k}_2)$ from (4.6) and
the conjectures.  Letting $k \rightarrow \infty$ we obtain $u(q_1)-\wt
u(q_1) = 0$.   \hfill  Q.E.D.

\bigskip
\underbar{Proof of Theorem 4.1}:  We begin by proving $u = \wt u$ under
the assumption that $u$ is even and $R < \lambda^D_1$.  Since $\Delta^ju$
satisfies the same hypotheses as $u$, we conclude from Lemma 4.4 that
$\Delta^j(u-\wt u)$ vanishes at all three boundary points, for any $j$.
Let $G(x,y)$ denote the Green's function and $G^j(x,y)$ the $j$-fold
iteration of $G$.  The vanishing at boundary points means that
$$
u(x) - \wt u(x) = \int G^j(x,y)\Delta^j(u(y)-\wt u(y))d\mu(y).
     \tag4.7
$$
We have an explicit representation
$$
G^j(x,y) = \sum^\infty_{k=1} (\lambda^D_k)^{-j}\varphi_k(x)\varphi_k(y)
     \tag4.8
$$
for an orthonormal basis of Dirichlet eigenfunctions $\{\varphi_k\}$ with
$-\Delta\varphi_k = \lambda^D_k\varphi_k$.  This yields the estimate
$$
\Big( \iint |G^j(x,y)|^2d\mu(x)d\mu(y)\Big)^{1/2}
     = \Big(\sum^\infty_{k=1} (\lambda^D_k)^{-2j}\Big)^{1/2}
      \le c(\lambda^D_1)^{-j}
     \tag4.9
$$
by the Weyl asymptotics of $\{\lambda^D_k\}$.  Thus
$$
\|u-\wt u\|_2
     \le c(\lambda^D_1)^{-j} \|\Delta^j(u-\wt u)\|_2
    \le c(\lambda^D_1)^{-j} R^j.
$$
Letting $j \rightarrow \infty$ we obtain $\|u-\wt u\|_2 = 0$ hence $u =
\wt u$ as desired.

Next we can remove the assumption that $u$ be even by writing $u$ as a
sum of even functions about each of the three boundary points using
(2.26).  It is clear that the hypotheses on $u$ are inherited by the
three summands, and a sum of three entire analytic functions is entire
analytic.

Finally, we need to relax the assumption that $R < \lambda^D_1$ to $R <
\lambda_2$.  To do this we consider $u \circ F_w$ for all words of length
2 (because $5^{-2}\lambda_2 < \lambda^D_2$).  Then $u \circ F_w$
satisfies (3.16) with $R < \lambda^D_1$, so by the previous argument it
is entire analytic.  This means for each $w$ there exists $\wt u_w$
entire analytic with $u = \wt u_w$ on $F_w(SG)$.  Next we claim that $\wt
u_{00} = \wt u_{01} = \wt u_{02}$.  To see this we may assume without
loss of generality that $\wt u_{00} = 0$ by replacing $u$ by $u-\wt
u_{00}$.  So $u$ is assumed to vanish on $F^2_0(SG)$, and we need to show
that it vanishes on $F_0(SG)$.  By Lemma 4.4 we have $u(F_0q_1) =
u(F_0q_2) = 0$, and more generally $\Delta^ju(F_0q_1) = \Delta^ju(F_0q_2)
= 0$ by the same reasoning for $\Delta^ju$.  Let us consider $\wt u_{01}$
which equals $u$ on $F_0F_1(SG)$.  At the point $F^2_0q_1$ where the
cells $F_0F_1(SG)$ and $F^2_0(SG)$ intersect, we have $\Delta^ju$
vanishing and also $\partial_n\Delta^ju$ vanishing (obvious for the
normal derivative with respect to $F^2_0(SG)$, and then true with respect
to $F_0F_1(SG)$ by the matching condition for normal derivatives).  Thus
the local power series expansion in $F_0F_1(SG)$ of $\wt u_{01}$ about
the point $F^2_0q_1$ contains only $P_{j3}$ terms, so $\wt u_{01}$ and
more generally $\Delta^j\wt u_{01}$ must be odd, so the vanishing of
$\Delta^j\wt u_{01}$ at the second boundary point $F_0q_1$ implies the
vanishing at the third boundary point $F_0F_1q_2$.  So our previous
argument shows that $\wt u_{01}$ is identically zero.

The same argument works in the other two cells of level one, so we now
know that there exist entire analytic functions $\wt u_0$, $\wt u_1$,
$\wt u_2$ such that $u = \wt u_j$ on $F_j(SG)$.  We need to show $\wt u_0
= \wt u_1 = \wt u_2$, and by subtracting $\wt u_0$ we may assume without
loss of generality that $\wt u_0 = 0$.  At this point we cannot simply
repeat the argument of the previous paragraph because the cell $F_1(SG)$
is too big.  Of course we can argue as before that $\wt u_1$ and more
generally $\Delta^j\wt u_1$ vanishes on all three boundary points
of $F_1(SG)$, and that it is odd about the vertex $F_0q_1$.  It is this
oddness that saves the argument.  Instead of (4.7) for $\wt u_1 \circ F_1$
we have
$$
\wt u_1 \circ F_1(x) = \int \wt G^j(x,y)\Delta^j(\wt u_1\circ F_1)(y)
          d\mu(y)
     \tag4.10
$$
where $\wt G^j$ denotes the $j$-fold iteration of the odd part of the
Green's function.  Instead of (4.8), $\wt G^j$ has the same
representation where the sum is restricted to the odd eigenfunctions.
The eigenfunction associated to $\lambda^D_1$ is even, so the smallest
eigenvalue appearing is $\lambda^D_2 \approx 55.8858.\ldots$ \ .  Thus we
obtain the estimate
$$
\|\wt u_1\circ F_1\|_2 \le c(\lambda^D_2)^{-j}5^{-j}R^j,
$$
and this shows $\wt u_1 = 0$ because $\lambda_2 \le 5\lambda^D_2$.
\hfill  Q.E.D.

\bigskip
It is interesting that the growth conditions (3.16) imply the specific
identities (2.22). There is nothing analogous to this in the theory of
real analytic functions.  In some way it is reminiscent of the Cauchy
integral formula for complex analytic functions.  But we don't want to
read too much into this, since (2.22) holds for nonanalytic functions as
well.

\proclaim{Corollary 4.5}  If $u$ is
defined on a cell $F_w(SG)$ and satisfies
$$
\|\Delta^ju\|_{L^\infty(F_w(SG))} = O(R^j) \ \text{for some} \
           R < \lambda_2
     \tag4.11
$$
then $u$ has a unique extension to an entire analytic function.
\endproclaim

\underbar{Proof}:  The theorem shows $u\circ F_w$ is entire analytic.
Then apply Theorem~3.6.~\hfill Q.E.D.

\bigskip
We can also consider entire analytic functions on any infinite blow-up of
SG.  The coefficients must satisfy (3.11) for all $R > 0$, and the
characterization requires the estimate (3.16) to hold locally for all $R >
0$.


\subheading{\S5. Expansions about junction points}

\medskip
A junction point is a boundary point of two cells, so an entire analytic
function will have two different local power series (3.14) centered at the
point, each valid in a different cell.  Since each local power series
determines the function, it also determines the other local power
series.  Since the coefficients of the local power series are just the
jets at the point with respect to each cell, these jets determine each
other.  The first goal of this seciton is to make this determination
explicit.

To be specific, consider the junction point $F_0q_1 = F_1q_0$.  We will
write $F_0q_1 = q_{01}$ and write
$$
(\Delta^ju(q_{01}), \partial_n\Delta^ju(q_{01}),
      \partial_T\Delta^ju(q_{01}))
   \tag5.1
$$
for the jet associated with the cell $F_0(SG)$, and $F_1q_0 = q_{10}$ and
$$
(\Delta^ju(q_{10}), \partial_n\Delta^ju(q_{10}),
      \partial_T\Delta^ju(q_{10}))
     \tag5.2
$$
for the jet associated with the cell $F_1(SG)$.  We know some
relationships between the jets (5.1) and (5.2), namely
$$
\Delta^ju(q_{01}) = \Delta^ju(q_{10}) \ \text{and} \
      \partial_n\Delta^ju(q_{01}) = -\partial_n\Delta^ju(q_{10}).
    \tag5.3
$$
Note that (5.3) is valid for all $u \in \dom \,\Delta^\infty$, but there
should be no connections between tangential derivatives without the
assumption that $u$ is an entire analytic function.  On the other hand,
for entire analytic functions, we expect an identity of the form
$$
\partial_T u(q_{01}) + \partial_Tu(q_{10})
     = \sum^\infty_{\ell=0} Y_\ell\partial_n\Delta^\ell u(q_{01})
     \tag5.4
$$
to hold for certain coefficients $Y_\ell$.  Note that (5.4) applied to
$\Delta^ju$ yields
$$
\partial_T\Delta^ju(q_{01}) + \partial_T\Delta^ju(q_{10})
     = \sum^\infty_{\ell=j} Y_{\ell-j} \partial_n\Delta^\ell u(q_{01}),
    \tag5.5
$$
and (5.3) and (5.5) show how the jets (5.1) and (5.2) determine each
other.  We may also interpret (5.4) as a matching condition for
tangential derivatives.

Our strategy for determining the $Y$ coefficients will be to first
consider the case when $u$ is a polynomial, making the sum finite.  It is
convenient to consider the monomials $P_{jk}^{(2)}$, because the $\rho_2$
symmetry is also a symmetry about $q_{01}$.  For even functions, both
sides of (5.4) are zero regardless of the $Y$ coefficients: the left side
vanishes because of the oddness of the tangential derivative, and the
right side because of the matching condition $\partial_n\Delta^\ell
u(q_{01}) = -\partial_n\Delta^\ell u(q_{01})$ and the evenness of the
normal derivative and Laplacian.  Thus we need only check (5.4) for the
monomials $P^{(2)}_{j3}$.

\proclaim{Lemma 5.1}  The matching condition (5.4) holds for all
polynomials for the $Y$ coefficients satisfying $Y_0 = 4$ and recursively
$$
\split
Y_j = -\alpha_j&-18 \sum^j_{\ell=0} n_{j+1-\ell}
          \frac{t_\ell}{5^\ell}
    + \sum^{j-1}_{\ell=0} Y_\ell\Big(\Big(\frac 32
          - \frac{5^{\ell-j}}{2}\Big)
           n_{j-\ell+1} \\
   &+ \sum^{j-\ell}_{k=0} (5\alpha_{j+1-\ell-k})n_k5^{-k}
        - 3n_{j+1-\ell-k} \alpha_k5^{-k})\Big) \ \text{for} \ j \ge 1.
\endsplit     \kern-1.5em\tag5.6
$$
\endproclaim

\underbar{Proof}:  When $j = 0$ we compute directly that
$\partial_TP^{(2)}_{03}(q_{01}) + \partial_TP^{(2)}_{03}(q_{10}) = -4$
and $\partial_nP^{(2)}_{03}(q_{01}) = -1$, so $Y_0 = 4$.  For $j \ge 1$
we use Corollary 2.14 to rearrange $P^{(2)}_{j3}$ around $q_0$. By (2.43)
we obtain
$$
P^{(2)}_{j3} = -\frac 12 P^{(0)}_{j3}
         + 3 \sum^j_{\ell=0} (\alpha_{j+1-\ell} P^{(0)}_{\ell 1}
              + n_{j+1-\ell} P^{(0)}_{\ell 2}).
    \tag5.7
$$
Because $P^{(2)}_{j3}$ is odd we have
$$
\partial_TP^{(2)}_{j3}(q_{01}) + \partial_TP^{(2)}_{j3}(q_{10})
      = 2\partial_TP^{(2)}_{j3}(q_{01}).
$$
By (5.7) and Theorem 2.12 we have
$$
2 \partial_TP^{(2)}_{j3}(q_{01})
     =  \alpha_j +18 \sum^j_{\ell=0} n_{j+1-\ell}
           \frac{t_\ell}{5^\ell}
     \tag5.8
$$
and
$$
\partial_n P^{(2)}_{j3}(q_{01})
     = \Big(\frac 32 - \frac 12 5^{-j}\Big)n_{j+1}
       + \sum^j_{k=0} (5\alpha_{j+1-k}n_k5^{-k}
        - 3n_{j+1-k}\alpha_k5^{-k}).
     \tag5.9
$$
Since $\Delta^\ell P^{(2)}_{j3} = P^{(2)}_{(j-\ell)3}$, we have that
(5.4) for $u = P^{(2)}_{j3}$ yields
$$
Y_j = \sum^{j-1}_{\ell=0} Y_\ell \partial_n P^{(2)}_{(j-\ell)3}(q_{01})
      - 2\partial_TP^{(2)}_{j3}(q_{01}).
$$
Substituting (5.8) and (5.9) yields (5.6).  \hfill  Q.E.D.

\bigskip
\proclaim{Conjecture 5.2}  The coefficients $Y_j$ satisfy
$$
|Y_j| \le c\lambda_2^{-j}.   \tag5.10
$$
\endproclaim

The numerical evidence for Conjecture 5.2 is presented in Table 5.1.


$$
\matrix
j    &Y_j      &(-\lambda_2)^jY_j \\
0    &  -4.   &-4. \\
1   &-0.08888888889   &12.05085573 \\
2   &0.0002304526749 & 4.235674447 \\
3   &-0.1434871749\  10^{-5}    &3.575397353  \\
4  & 0.1023938272\  10^{-7}     &3.459038654  \\
5   &-0.7503519662\  10^{-10}   &3.436505741  \\
6  & 0.5527533783 \ 10^{-12}    &3.432052039  \\
7   &-0.4076138308 \ 10^{-14}   &3.431166398  \\
8  &0.3006465014 \ 10^{-16}     &3.430989845  \\
9  &-0.2217590148 \ 10^{-18}    &3.430954602  \\
10  &0.1635723837 \ 10^{-20}    &3.430947563  \\
11  &-0.1206533528 \ 10^{-22}   &3.430946155  \\
12  &0.8899568485\  10^{-25}    &3.430945874  \\
13  &-0.6564452839 \ 10^{-27}   &3.430945818  \\
14  &0.4842037197 \ 10^{-29}    &3.430945807  \\
15  &-0.3571558034\  10^{-31}   &3.430945805  \\
16  &0.2634433871 \ 10^{-33}    &3.430945805  \\
17  &-0.1943197270\  10^{-35}   &3.430945805  \\
18  &0.1433330961\  10^{-37}    &3.430945805  \\
19  &-0.1057246052 \ 10^{-39}   &3.430945805  \\
20  &0.7798402782 \ 10^{-42}   &3.430945805
\endmatrix
$$

\centerline{Table 5.1}


\proclaim{Theorem 5.3}  Assume Conjecture 5.2.  If $u$ is
any entire analytic function, then (5.4) and (5.5) hold for the $Y$
coefficients given in Lemma 5.1.  More generally, if $x$ is any junction
point in $V_{m+1}\setminus V_m$, then
$$
\partial_T\Delta^ju(x) + \partial_T^*\Delta^ju(x)
    = \sum^\infty_{\ell=j} 3^m5^{-m(\ell-j)} Y_{\ell-j}
       \partial_n\Delta^\ell u(x),
    \tag5.11
$$
where $\partial_T$ and $\partial_n$ are derivatives with respect to the
left cell at $x$ and $\partial_T^*$ is the derivative with respect to the
right cell.
\endproclaim

\underbar{Proof}:  Note that the right side of (5.4) converges
absolutely.  The issue is then whether the term--by--term differentiation
of power series extends to normal and tangential derivatives at points
other than the expansion point.  For normal derivatives this is easy to
see because of the integral representation.  But in any case this follows
by combining Theorem 3.4 (the explicit expression (3.12) for the
rearranged coefficients) with Theorem 3.1 (the jet formula (3.5) at the
expansion point). We then obtain (3.10) by applying (3.5) to the function
$u\circ F_w$ for $|w| = m$.  \hfill  Q.E.D.

\bigskip
Next we consider the question of what would be a natural notion of a
power series expansion centered about a junction point.  We will see that
there is no completely satisfactory answer.  Again to be specific we
consider the point $q_{01} = q_{10}$.  We would like to have at least the
following four conditions holding:

(i) every entire analytic function has an expansion;

(ii) the expansion is valid in a neighborhood of $q_{01}$, perhaps
$F_0(SG) \cup F_1(SG)$;

(iii) the individual terms are polynomials that vanish to higher and
higher order near $q_{01}$;

(iv)  the rate of growth of the coefficients should be characterized for
entire analytic functions.

The local power series with respect to one of the cells, say $F_0(SG)$,
gives a satisfactory answer only on that cell, but if we continue those
monomials around we will find that the vanishing rate near $q_{10}$ is
not satisfactory.  In fact the tangential derivatives will have to be
nonzero by Lemma 5.1.  For this reason we consider carefully what it
takes to meet condition (iii). We denote by $P^{(01)}_{jk}$ the monomials
of the $F_0(SG)$ local power series about $q_{01}$, so that
$$
\split
&\Delta^\ell P^{(01)}_{jk}(q_{01}) = \delta_{j\ell}\delta_{k1} \\
&\partial_n\Delta^\ell P^{(01)}_{jk}(q_{01}) = \delta_{j\ell}\delta_{k2}\\
&\partial_T\Delta^\ell P^{(01)}_{jk}(q_{01}) = \delta_{j\ell}\delta_{k3}
\endsplit
$$
or more precisely
$$
\split
&P^{(01)}_{j1}(x) = 5^{-j}P^{(1)}_{j1}(F^{-1}_0x) \\
&P^{(01)}_{j2}(x) = \frac 35 5^{-j}P^{(1)}_{j2}(F_0^{-1}x) \\
&P^{(01)}_{j3}(x) = 5^{-j-1}P^{(1)}_{j3}(F_0^{-1}x).
\endsplit
$$
Note that $P^{(01)}_{j1}$ and $P^{(01)}_{j3}$ extend to even polynomials
about $q_{01}$, so they will have the same vanishing rate on both cells.
We want to replace $P^{(01)}_{j2}$ by a different polynomial $\wt
P^{(01)}_{j2}$ that will have the same $j$--jet (except for
$\partial_T\Delta^ju(q_{01})$), but will extend to be odd.  This will
give it the correct order of vanishing, but in exchange we have to take
a higher order polynomial.  The lowest possible order is $2j$:
$$
\wt P^{(01)}_{j2} = \sum^j_{\ell=0} (a_{j(j-\ell)} P^{(01)}_{(j+\ell)2}
         + b_{j(j-\ell)} P^{(01)}_{(j+\ell)3})
   \tag5.12
$$
for the appropriate choice of constants.  Note that we can exclude
$P^{(01)}_{(j+\ell)1}$ terms because we want the possibility of odd
extension.  We will take $a_{jj} = 1$ in order to obtain the correct
$j$--jet.  The odd extension means $\partial_T\Delta^n\wt
P^{(01)}_{j2}(q_{01}) = \partial_T\Delta^n\wt P^{(01)}_{j2}(q_{10})$, so
we have $2j+1$ equations of the form (5.5) to satisfy, and these
will determine the remaining $2j+1$ constants.  The equations are
$$
2 \partial_T \Delta^n\wt P^{(01)}_{j2}(q_{10})
    = \sum^{2j}_{k=n} Y_{k-n}\partial_n\Delta^k \wt P^{(01)}_{j2}(q_{02}),
   \tag5.13
$$
and when $0 \le n < j$ the left side is zero and we obtain
$$
0 = \sum^{2j}_{k=n} Y_{k-n} \partial_n \Delta^k \wt P_{j2}^{(01)}(q_{01})
       = \sum^{2j}_{k=j} Y_{k-n} a_{j(2j-k)}
$$
so
$$
0 = \sum^j_{\ell =0} Y_{2j-\ell-n} a_{j\ell}.
   \tag5.14
$$
We use these equations to solve for $a_{j\ell}$.  When $n \le j \le 2j$
the left side of (5.13) is $2 b_{j(2j-n)}$ so
$$
2b_{j(2j-n)} = \sum^{2j}_{k=n} Y_{k-n} a_{j(2j-k)},
$$
and by letting $\ell = 2j-n$ we have
$$
b_{j\ell} = \frac 12 \sum^\ell_{k=0}Y_{k} a_{j(\ell-k)} \ \text{for} \
       0 \le \ell \le j.
   \tag5.15
$$

In Table 5.2 we show the values of $a_{j\ell}$ and $b_{j\ell}$ for small
values of $j$.  It is difficult to discern a pattern in these results.  We
have obtained graphs of $\wt P_{j2}^{(01)}$ for small values of $j$ using
(5.12), but it appears that round--off error becomes significant before
any pattern emerges, so we are not able to offer any conjectures about the
growth rate of these functions as $j \rightarrow \infty$.

\newpage
\def\s{\eightpoint}
$$
\matrix
\s j   &\s l   &\s a_{jl}    &\s b_{jl} &   &\s j  &\s l  &\s a_{jl}
              &\s b_{jl} \\
\s 0   &\s 0   &\s 1.   &\s 2.   &
        &\s 7   &\s 0   &\s 0.1330959781\  10^{23}
               &\s 0.2661919562 \ 10^{23} \\
\s 1   &\s 0   &\s 0.02252966406   &\s 0.04505932812  &
       &\s 7   &\s 1   &\s 0.6141913960\ 10^{21}
                &\s 0.7847295317\ 10^{21}\\
\s 1   &\s 1   &\s 1.   &\s 1.999249011   &
        &\s 7  &\s 2  &\s 0.6084736857\ 10^{19}
             &\s 0.1968365718\ 10^{23}\\
\s 2   &\s 0   &\s 6461.417615   &\s 12922.83523   &
     &\s 7  &\s 3  &\s 0.2707503937 \ 10^{17}
              &\s 0.1030137030 \ 10^{22}\\
\s 2   &\s 1   &\s -39.86777272   &\s -295.1161326 &
       &\s 7  &\s 4  &\s 0.4581523610\ 10^{14}
              &\s 0.1231428577 \ 10^{20}\\
\s 2   &\s 2   &\s 1.   &\s 9563.195714   &
      &\s 7  &\s 5     &\s 0.2620127789 \ 10^{11}
             &\s 0.5414059059 \ 10^{17}\\
\s 3   &\s 0   &\s 0.1631072895 \ 10^{7}  &\s 0.3262145790 \ 10^{7} &
       &\s 7  &\s 6  &\s 3880.162356   &\s 0.9158266227 \ 10^{14}  \\
\s 3   &\s 1   &\s 48581.69671    &\s 42794.29693   &
         &\s 7  &\s7  &\s 1.  &\s 0.5243561927 \ 10^{11}  \\
\s 3   &\s 2   &\s -109.6002902   &\s 0.2411384099\  10^{7}   &
      &\s 8  &\s 0  &\s -0.2849367688 \ 10^{25}
            &\s -0.5698735375 \ 10^{25}\\
\s 3   &\s 3   &\s 1.   &\s 86782.07999    &
       &\s 8  &\s 1  &\s -0.1352864496\ 10^{24}
         &\s 0.1755939762 \ 10^{24}  \\
\s 4   &\s 0   &\s -0.1623039023 \ 10^{10}   &\s -0.3246078045\  10^{10}  &
      &\s 8  &\s 2  &\s -0.1478090302 \ 10^{22}
        &\s -0.4214174069 \ 10^{25} \\
\s 4   &\s 1   &\s -0.6442287860 \ 10^{8}  &\s -0.7474445645\  10^{8}  &
       &\s 8  &\s 3  &\s -0.7540725789 \ 10^{19}
            &\s -0.2261525660 \ 10^{24}\\
\s 4   &\s 2   &\s -299734.8354   &\s -0.2399788368 \ 10^{10}  &
       &\s 8  &\s 4  &\s -0.1760661536 \ 10^{17}
       &\s -0.2930516976 \ 10^{22} \\
\s 4   &\s 3   &\s -347.4611669   &\s -0.1101312661 \ 10^{9} &
     &\s 8  &\s 5  &\s -0.1895987908 \ 10^{14}
         &\s -0.1511819510 \ 10^{20} \\
\s 4   &\s 4   &\s 1.   &\s -751724.7199   &\
          &\s 8  &\s 6  &\s -0.7756675150 \ 10^{10}
          &\s -0.3520528934 \ 10^{17}\\
\s 5   &\s 0   &\s 0.1010368178 \ 10^{14}  &\s 0.2020736356 \ 10^{14}
       &  &\s 8  &\s 7  &\s -3618.462380
             &\s -0.3790729379 \ 10^{14}  \\
\s 5   &\s 1   &\s 0.4380632964\  10^{12}  &\s 0.5393372002\  10^{12}
      &  &\s 8  &\s 8  &\s 1.  &\s -0.1552676258 \ 10^{11} \\
\s 5   &\s 2   &\s 0.3374174349\  10^{10}  &\s 0.1494085527 \ 10^{14}
     &  &\s 9  &\s 0  &\s 0.4817483229 \ 10^{29}
        &\s 0.9634966458 \ 10^{29}\\
\s 5   &\s 3   &\s 0.1015644445\  10^{8}   &\s 0.7403235769 \ 10^{12}
      &  &\s 9  &\s 1  &\s 0.2289760048 \ 10^{28}
         &\s 0.2973692352 \ 10^{28} \\
\s 5   &\s 4   &\s -909.3198857   &\s 0.7249040413 \ 10^{10}   &
       &\s 9  &\s 2  &\s 0.2513117964 \ 10^{26}
           &\s 0.7125008828 \ 10^{29}\\
\s 5   &\s 5   &\s 1.   &\s 0.1921254540 \ 10^{8} &
      &\s 9  &\s 3  &\s 0.1299020030 \ 10^{24}
            &\s 0.3827224251 \ 10^{28}\\
\s 6   &\s 0   &\s -0.1389829261 \ 10^{18}  &\s -0.2779658521 \ 10^{18}
     &  &\s 9  &\s 4  &\s 0.3154544064 \ 10^{21}
         &\s 0.4977745865 \ 10^{26}\\
\s 6   &\s 1   &\s -0.6247328496\  10^{16}  &\s -0.7861892790 \ 10^{16}
       &  &\s 9  &\s 5  &\s 0.3859718201 \ 10^{18}
           &\s 0.2600348690 \ 10^{24} \\
\s 6   &\s 2   &\s -0.5605362673\  10^{14}  &\s -0.2055333917 \ 10^{18}
     &  &\s 9  &\s 6  &\s 0.2325380299 \ 10^{15}
       &\s 0.6310751388 \ 10^{21}  \\
\s 6   &\s 3   &\s -0.2151475440\  10^{12}  &\s -0.1051115464\  10^{17}
     &  &\s 9  &\s 7  &\s 0.5539946952 \ 10^{11}
        &\s 0.7717084596 \ 10^{18} \\
\s 6   &\s 4   &\s -0.2169919676\  10^{9}  &\s -0.1159983908 \ 10^{15}
     &  &\s 9  &\s 8  &\s -6592.977986
            &\s 0.4652032965\ 10^{15} \\
\s 6   &\s 5   &\s -1787.130925   &\s -0.4257054009\  10^{12}  &
       &\s 9  &\s 9  &\s 1.  &\s 0.1107495241 \ 10^{12} \\
\s 6   &\s 6   &\s 1.   &\s -0.4383706038\  10^{9}
\endmatrix
$$

\bigskip
\centerline{Table 5.2}
\newpage

\subheading{\S6.  Exponentials}

\medskip
Eigenfunctions of the Laplacian give us a natural class of special
functions on SG.  Until now, most attention has been paid to
eigenfunctions satisfying Dirichlet or Neumann boundary conditions, which
forces the eigenvalue to be positive.  In contrast, we will mainly
explore negative eigenvalues in this section, so we are exploring the
analog of the functions $\cosh\sqrt\lambda t$ and $\sinh\sqrt\lambda t$
on the unit interval and their extension to the positive real line.  Of
particular interest is the linear combination that yields
$e^{-\sqrt\lambda t}$, the unique choice that exhibits exponential decay
(either as $\lambda \rightarrow \infty$ or as $t \rightarrow \infty$) as
opposed to exponential growth.  It is embarrassing to note that the
exponential $e^{\sqrt\lambda t}$ does not distinguish itself among linear
combinations of $\cosh\sqrt\lambda t$ and $\sinh\sqrt\lambda t$, if one
is forbidden to use odd order derivatives.  So we have not been able to
find its analog on SG.

The space of all eigenfunctions with a fixed eigenvalue has dimension
three, as long  as one avoids Dirichlet eigenvalues.  For fixed $\lambda
>  0$ we can choose a basis $C_\lambda$, $S_\lambda$, $Q_\lambda$ for the
space of solutions to
$$
-\Delta u = -\lambda u
    \tag6.1
$$
determined by the conditions that $C_\lambda$ and $S_\lambda$ are even
and $Q_\lambda$ is odd with respect to $\rho_0$, and
$$
C_\lambda(q_0) = 1, \quad \partial_nC_\lambda(q_0) = 0   \tag6.2
$$
$$
S_\lambda(q_0) = 0, \quad \partial_nS_\lambda(q_0) = a_\lambda   \tag6.3
$$
$$
\partial_TQ_\lambda(q_0) = 1    \tag6.4
$$
where the normalization factor $a_\lambda$ will be chosen later.  This
means that we have global power series representation
$$
C_\lambda(x) = \sum^\infty_{j=0} \lambda^j P^{(0)}_{j1}(x)
    \tag6.5
$$
and
$$
Q_\lambda(x) = \sum^\infty_{j=0} \lambda^j P^{(0)}_{j3}(x),
    \tag6.6
$$
and a local power series representation
$$
S_\lambda(x) = a_\lambda \sum^\infty_{j=0} \lambda^j P^{(0)}_{j2}(x)
    \tag6.7
$$
valid on $F^n_0(SG)$ provided $\lambda < 5^n\lambda_2$.  We may also use
(6.5) and (6.6) on the blowups $F_0^{-n}(SG)$ for any $n$.  Of course,
none of these functions are entire analytic for $\lambda \ge \lambda_2$.

We will consider the infinite blowup $SG_\infty = \ds\bigcup^\infty_{n=0}
F_0^{-n}(SG)$ to play the role of the positive reals vis-a-vis the unit
interval.  Of course there are uncountably many infinite blow-ups of SG.
We have chosen the simplest one to study first.  To understand the
``behavior at infinity'' of these functions it suffices to study the
values at the points $x_n = F^n_0q_1$ as $n \rightarrow -\infty$, for we
may then get the values at the points $y_n = F^n_0q_2$ by parity, and
then fill in by spectral decimation.

For $SG_\infty$ we have graphs $\Gamma_n$ for any integer $n$.  Since
$-\lambda$ is negative we never encounter the exceptional eigenvalues 2,
5 and 6.  Thus the method of spectral decimation says that $u$ satisfies
(6.1) on $SG_\infty$ if and only if the restriction of $u$ to $\Gamma_n$
is a graph eigenfunction with eigenvalue $\lambda_n$, where
$\{\lambda_n\}_{n\in\Bbb Z}$ is a sequence of negative numbers
characterized by
$$
\lambda_{n-1} = \lambda_n(5-\lambda_n)    \tag6.8
$$
and
$$
-\lambda = \lim_{n\rightarrow\infty} \frac 32 5^n\lambda_n.   \tag6.9
$$
Note that $\lambda_n \rightarrow 0$ as $n \rightarrow \infty$ and
$\lambda_n \rightarrow -\infty$ as $n \rightarrow -\infty$.  It is easy
to see that the sequence $\{\lambda_j\}$ is uniquely characterized by
these conditions, and the values may be effectively computed to any
desired accuracy by replacing the limit in (6.9) by the value for a fixed
large $n$ and then using (6.8) to run $n$ down.

The fact that $u$ restricted to $\Gamma_n$ is a $\lambda_n$-eigenfunction
means that if we take any cell of level $n-1$ with boundary points $a$,
$b$, $c$, and if $d$ is the midpoint between $a$ and $b$, then
$$
u(d) = \frac{(4-\lambda_n)(u(a)+u(b))+2u(c)}
           {(2-\lambda_n)(5-\lambda_n)}
     \tag6.10
$$
(see [DSV] Algorithm 2.4).

\proclaim{Lemma 6.1}  The recurrence relations
$$
C_\lambda(x_n)
    = \frac{(4-\lambda_n)+(6-\lambda_n)C_\lambda(x_{n-1})}
           {(2-\lambda_n)(5-\lambda_n)}
     \tag6.11
$$
$$
S_\lambda(x_n) = \frac{(6-\lambda_n)S_\lambda(x_{n-1})}
           {(2-\lambda_n)(5-\lambda_n)}
    \tag6.12
$$
and
$$
Q_\lambda(x_n) = \frac{Q_\lambda(x_{n-1})}{5-\lambda_n}
    \tag6.13
$$
hold for all integers $n$.
\endproclaim

\underbar{Proof}:  Apply (6.10) for $a = q_0$, $b = F_0^{n-1}(q_1)$,
$c = F_0^{n-1}(q_2)$ and $d = F_0^n(q_1)$.~\hfill  Q.E.D.

\bigskip
\proclaim{Lemma 6.2}  The function $C_\lambda$ is positive.  The function
$S_\lambda$, with the appropriate choice of $a_\lambda$, is positive
everywhere except at $q_0$ where it vanishes.  The function $Q_\lambda$
vanishes on the symmetry line through $q_0$ and is positive on the $q_1$
half of the symmetry line.
\endproclaim

\underbar{Proof}:  Because $\lambda_n < 0$ for all $n$, the coefficients
in (6.10-6.13) are all positive.  That means that if $u$ is nonnegative on
the boundary of a cell and strictly positive at one of the boundary points
then it is strictly positive in the interior.  Thus it suffices to show
that $C_\lambda(x_n)$, $S_\lambda(x_n)$ and $Q_\lambda(x_n)$ are
positive.  For $S_\lambda$ and $Q_\lambda$ it suffices to show
$S_\lambda(q_1)$ and $Q_\lambda(q_1)$ are positive, since we can solve
(6.12) and (6.13) for $S_\lambda(x_{n-1})$ and $Q_\lambda(x_{n-1})$ with
positive coefficients.  But we can make $S_\lambda(q_1) > 0$ by the
appropriate choice of sign (negative) for $a_\lambda$, and
$Q_\lambda(q_1) > 0$ follows easily from $\partial_TQ_\lambda(q_0) = 1$.
When we solve (6.11) we obtain
$$
C_\lambda(x_{n-1})=
\frac{(2-\lambda_n)(5-\lambda_n)C_\lambda(x_n)-(4-\lambda_n)}
          {6-\lambda_n}\ ,
    \tag6.14
$$
which contains a negative coefficient.  Nevertheless, if $C_\lambda(x_n)
>  1$ then (6.14) implies
$$
C_\lambda(x_{n-1}) > \frac{(2-\lambda_n)(5-\lambda_n)-(4-\lambda_n)}
        {6-\lambda_n} > 1,
$$
so it suffices to show $C_\lambda(q_1) > 1$.  This follows because the
contrary assumption $C_\lambda(q_1) \le 1$ and (6.13) would imply
$\partial_nC_\lambda(q_0) > 0$.   \hfill  Q.E.D.

\bigskip
\proclaim{Theorem 6.3}  (a) For all $n$ we have
$$
C_\lambda(x_n) = 1 - \frac{\lambda_n}{4}\ .     \tag6.15
$$

(b) For the appropriate choice of $a_\lambda$ we have
$$
S_\lambda(x_n) = -\frac{\lambda_n}{4} \prod^\infty_{k=0}
         \Big(1 + \frac{4}{2-\lambda_{n-k}}\Big),
      \tag6.16
$$
and hence
$$
\lim_{n\rightarrow-\infty} S_\lambda(x_n)/C_\lambda(x_n) = 1.
       \tag6.17
$$

(c) For all $n < 0$ we have
$$
Q_\lambda(x_n) = -\frac 34 \ \frac{\lambda_n}{\lambda}   \tag6.18
$$
and hence
$$
\lim_{n\rightarrow-\infty} Q_\lambda(x_n)/C_\lambda(x_n)
      = \frac 3\lambda.    \tag6.19
$$
\endproclaim

\underbar{Proof}:  (a) A direct calculation using (6.8) shows that
$1-\frac{\lambda_n}{4}$ satisfies the same recurrence relation (6.11) as
$C_\lambda(x_n)$.  Thus if we define $\wt C_\lambda(x_n) = 1 -
\frac{\lambda_n}{4}$, $\wt C_\lambda(q_0) = 1$ and extend $\wt C_\lambda$
to all of $SG_\infty$ using (6.10), we will have an even
$\lambda$-eigenfunction.  But a direct computation shows
$$
\partial_n\wt
C_\lambda(q_0) = \lim_{j\rightarrow\infty} \Big(\frac 53\Big)^j
         \frac 12 \lambda_j = 0
$$
because $\lambda_j = O(5^{-j})$ as $j \rightarrow \infty$.  So $\wt
C_\lambda = C_\lambda$, proving (6.15).

(b) First we observe that the infinite product in (6.16) converges,
because of the rapid growth of $\lambda_n$ as $n \rightarrow -\infty$.
Since (6.12) may be written (using (6.8))
$$
\frac{S_\lambda(x_n)}{\lambda_n}
     = \Big(1 + \frac{4}{2-\lambda_n}\Big)
           \frac{S_\lambda(x_{n-1})}{\lambda_{n-1}}\ ,
    \tag6.20
$$
it follows that the right side of (6.16) satisfies (6.12).  Since
$S_\lambda$ was only defined up to a multiplicative constant, we may
choose $a_\lambda$ to make (6.16) hold.  Note that from (6.20) we obtain
$S_\lambda(x_n) = O\big(\big(\frac 35\big)^n\big)$ as $n \rightarrow
\infty$, which is consistent with $S_\lambda(q_0) = 0$ and
$\partial_nS_\lambda(q_0) \neq 0$.  Then (6.17) follows from (6.15) and
(6.16) by inspection.

(c) We may rewrite (6.13) as
$$
\frac{Q_\lambda(x_n)}{\lambda_n}
    = \frac{Q_\lambda(x_{n-1})}{\lambda_{n-1}}
$$
using (6.8), hence $Q_\lambda(x_n) = \lambda_nQ_\lambda(x_0)$ for all
$n$.  But then
$$
\split
1 = \partial_TQ_\lambda(q_0)
    &= \lim_{n\rightarrow\infty} 5^n(Q(x_n)-Q_\lambda(y_n)) \\
   &= 2Q_\lambda(x_0)\lim_{n\rightarrow\infty} 5^n\lambda_n \\
   &= -\frac 43 \lambda Q_\lambda(x_0).
\endsplit
$$
This proves (6.18), and then (6.19) follows by inspection.  \hfill Q.E.D.

\bigskip
We can compute the value of $a_\lambda = \partial_nS_\lambda(q_0)$
exactly.  From the definition and (6.16) we have
$$
\split
\partial_nS_\lambda(q_0)
    &= -2\lim_{n\rightarrow\infty} \Big(\frac 53\Big)^n S_\lambda(x_n)\\
   &= \lim_{n\rightarrow\infty} \frac{\lambda_n}{2}
       \Big(\frac 53\Big)^n \prod^\infty_{k=0}
           \Big(1 + \frac{4}{2-\lambda_{n-k}}\Big) \\
   &= -\frac 13 \lambda \lim_{n\rightarrow\infty} \frac{1}{3^n}
         \prod^\infty_{k=0} \Big(1 + \frac{4}{2-\lambda_{n-k}}\Big) \\
   &= -\frac 13 \lambda \prod^\infty_{j=0}
        \Big(1 + \frac{4}{2-\lambda_{-j}}\Big)\lim_{n\rightarrow\infty}
           \prod^n_{k=1} \Big(\frac{6-\lambda_k}{6-3\lambda_k}\Big)\\
   &= -\frac 13 \lambda \prod^\infty_{j=0}
        \Big(1 + \frac{4}{2-\lambda_{-j}}\Big)\prod^\infty_{k=1}
         \Big(\frac{6-\lambda_k}{6-3\lambda_k}\Big).
\endsplit     \tag6.21
$$

\medskip
\noindent
\underbar{Definition 6.4}:  For $\lambda < 0$ define the {\it decaying
exponential} function $E_\lambda$ by
$$
E_\lambda(x) = C_\lambda(x) - S_\lambda(x).  \tag6.22
$$

\proclaim{Theorem 6.5}  $E_\lambda(x_n) = O(\lambda_n^{-1})$ as $n
\rightarrow -\infty$.  In fact
$$
\lim_{n\rightarrow -\infty} \lambda_nE_\lambda(x_n) = -1  \tag6.23
$$
and
$$
\lim_{n\rightarrow-\infty} C_\lambda(x_n)^2 - S_\lambda(x_n)^2 = \frac 12.
      \tag6.24
$$
More precisely
$$
E_\lambda(x_n) = \frac{2}{2-\lambda_n}
                  + \frac{\lambda_n}{2-\lambda_{n-1}}
           + \frac{4\lambda_n}{(2-\lambda_n)(2-\lambda_{n-1})}
     + O(\lambda_n^{-3}).
     \tag6.25
$$
\endproclaim

\underbar{Proof}:  From (6.16) we obtain
$$
S_\lambda(x_n) = -\frac{\lambda_n}{4}
        \Big(1 + \frac{4}{2-\lambda_n}\Big)
         \Big(1 + \frac{4}{2-\lambda_{n-1}}\Big)
     + O(\lambda_n^{-3})
     \tag6.26
$$
because $\lambda_n/\lambda_{n-2} = O(\lambda_n^{-3})$.  Substituting
(6.26) into (6.22) and using (6.15) we obtain (6.25).  Using (6.8) we see
that the first two terms on the right side of (6.25) sum to
$$
\frac{2}{2-\lambda_n} + \frac{\lambda_n}{2-5\lambda_n+\lambda^2_n}
      = -\frac{1}{\lambda_n} + O(\lambda_n^{-2}).
$$
The third term is clearly $O(\lambda_n^{-2})$, so we obtain (6.23).  From
(6.26) we find $S_\lambda(x_n) = -\frac{\lambda_n}{4} + O(1)$ and this
yields (6.24).  \hfill  Q.E.D.

\bigskip
Note that (6.26) and (6.25) allow for the efficient computation of
$S_\lambda$ and $E_\lambda$ for $n$ sufficiently negative.  On the other
hand (6.22) is computationally unstable since it involves subtracting
values that are large and nearly identical.  In Table 6.1 we present some
numerical computations of these functions.

Instead of fixing $\lambda$ and taking the limit as $n \rightarrow
-\infty$, we could look at values at $x_0$ and let $\lambda \rightarrow
-\infty$.  As long as $|\lambda_0|$ is large, (6.25) and (6.26) will be
good estimates.  Table 6.2 shows this behavior.  We could also allow
$\lambda$ to be complex, as long as the real part is positive to avoid
the exceptional values for $\lambda_n$.

We now turn our attention to eigenfunctions with positive eigenvalues,
with the goal of using information gleaned from spectral decimation to
shed some light on the recursion relations from Section 2.  Keeping the
same notation as before, we are interested in the function
$$
C_{-\lambda}(x) = \sum^\infty_{j=0} (-\lambda)^j P_{j1}(x)
$$
$$
\matrix
-j    &\lambda_{-j}   &C_\lambda(x_{-j})   &S_\lambda(x_{-j}) \\
                   \vspace{.07in}
  0     &-10.      &3.500000000        &3.421641174 \\
-1    &-150.     &38.50000000        &38.49346321 \\
-2    &-23250.   &5813.500000        &5813.499957  \\
-3    &-.540678750\  10^{9}   &.1351696885\  10^{9}
                    &.1351696885\  10^{9} \\
-4    &-.2923335134 \ 10^{18}    &.7308337835 \ 10^{17}
                &.7308337835\  10^{17} \\
-5   &-.8545888306\  10^{35}    &.2136472076 \ 10^{35}
               &.2136472076\  10^{35} \\
-6   &-.7303220694\  10^{70}    &.1825805173\  10^{70}
              &.1825805173 \ 10^{70} \\
-7   &-.5333703250\  10^{140}   &.1333425813\  10^{140}
             &.1333425813\  10^{140} \\
-8   &-.2844839036 \ 10^{280}    &.7112097590 \ 10^{279}
                &.7112097590 \ 10^{279}  \\
-9   &-.8093109142\  10^{559}   &.2023277285 \ 10^{559}
               &.2023277285 \ 10^{559} \\
-10   &-.6549841558\  10^{1118}   &.1637460389\  10^{1118}
               &.1637460389\  10^{1118} \\  \vspace{.1in}
-j   &Q_\lambda(x_{-j})   &E_\lambda(x_{-j})
                 &\lambda_{-j}E_\lambda(x_{-j}) \\   \vspace{.07in}
  0   &.7008295323     &.07835882554       &-.7835882554 \\
-1   &10.51244298     &.006536787301      &-.9805180952 \\
-2   &1629.428662     &.00004300520387    &-.9998709899 \\
-3   &.3789236353\  10^{8}   &.1849527089 \ 10^{-8}
                 &-.9999999945 \\
-4  &.2048759594 \ 10^{17}    &.3420750458\  10^{-17}
              &-1.0000000000 \\
-5  &.5989210902\  10^{34}    &.1170153370\ 10^{-34}
             &-1.0000000000 \\
-6  &.5118312741\  10^{69}    &.1369258909\  10^{-69}
             &-1.0000000000 \\
-7  &.3738016753\  10^{139}    &.1874869960 \ 10^{-139}
            &-1.0000000000 \\
-8  &.1993747210\  10^{279}    &.3515137367\  10^{-279}
            &-1.0000000000 \\
-9  &.5671889891\  10^{558}    &.1235619071\  10^{-558}
            &-1.0000000000 \\
-10 &.4590322393\ 10^{1117}    &.1526754489\  10^{-1117}
           &-1.0000000000 \\
\endmatrix
$$

\centerline{Table 6.1.  Values of functions at $x_{-j}$ for
$\lambda = 10.70160380$.}

\bigskip
\noindent
and its values at the special points $x_0 = q_1$ and $x_1 = F_0q_1$.  It
is convenient to define $\lambda_n$ (here we only care about $n \ge 0$)
to satisfy (6.8) but to remove the minus sign in (6.9).  For the
Dirichlet and Neumann eigenfunctions we know exactly what these values
are, and then we can use Theorem 6.3 (a) to conclude that
$C_{-\lambda}(x_0) = 1-\frac{\lambda_0}{4}$ and $C_{-\lambda}(x_1) = 1
-\frac{\lambda_1}{4}$.  (Strictly speaking, we need to use an analytic
continuation and limit argument to get this for the values we are
interested in.)  In particular, if $\lambda_0 = -6$ then
$C_{-\lambda}(x_0) = 5/2$, or
$$
\sum^\infty_{j=0} (-\lambda)^j P_{j1}(q_1)
     = \sum^\infty_{j=0} (-\lambda)^j \alpha_j = 5/2.
$$
$$
\matrix
\lambda_0     &\lambda      &E_\lambda(x_0)
           &\text{first 2 terms}
            &\text{first 3 terms} \\
   &    &  &\text{in (6.25)}  &\text{in (6.25)} \\ \vspace{.08in}
-100    &44.19536761    &.009711493217    &.01008584733
           &.009712435727 \\
-500    &87.71437197    &.001988095160    &.002003881410
           &.001988103065 \\
-1000   &112.0105482    &.0009970119472    &.001000985089
             &.0009970129413 \\
-5000   &182.0354932    &.0001998800959    &.0002000398801
            &.0001998801039 \\
-10000   &218.2833208    &.00009997001199      &.0001000099850
             &.00009997001299 \\
-50000   &317.2473555    &.00001999880010    &.00002000039988
            &.00001999880010 \\  \vspace{.2in}
\lambda_0     &\lambda      &S_\lambda(x_0)
           &\text{first 2 factors}
            &\text{first 3 factors} \\
   &    &  &\text{in (6.26)}  &\text{in (6.26)} \\ \vspace{.08in}
-100   &44.19536761    &25.99028851    &25.98039216    &25.99028756 \\
-500   &87.71437197    &125.9980119    &125.9960159    &125.9980119 \\
-1000   &112.0105482   &250.9990030    &250.9980040    &250.9990030 \\
-5000   &182.0354932   &1250.999800    &1250.999600    &1250.999800 \\
-10000  &218.2833208   &2500.999900   &2500.999800    &2500.999900 \\
-50000  &317.2473555   &12500.99998   &12500.99996    &12500.99998 \\
\endmatrix
$$

\centerline{Table 6.2. Values of functions at $x_0$ for various $\lambda$
values.}

\bigskip
\noindent
This happens when $\lambda = \lambda_2$, the second nonzero
Neumann eigenvalue (not to be confused with the $\lambda_2$ in (6.8) and
(6.9)).  This allows us to compute the limit of $\beta_j/t_{j+1}$ as $j
\rightarrow \infty$.  Indeed, from (2.34) we
have
$$
\frac{\beta_j}{t_{j+1}}
    = 6 \sum^j_{\ell=0} \alpha_{j+1-\ell}\Big(\frac{t_\ell}{t_{j+1}}\Big)
    = 6 \sum^{j+1}_{\ell=0} \alpha_\ell
      \Big(\frac{t_{j+1-\ell}}{t_{j+1}}\Big)-6 .
$$
We expect to have
$$
\frac{t_{j+1-\ell}}{t_{j+1}} \approx (-\lambda_2)^\ell
$$
and so
$$
\lim_{j\rightarrow\infty} \frac{\beta_j}{t_{j+1}}
    = 6 \sum^\infty_{\ell=0} \alpha_\ell(-\lambda_2)^\ell - 6
    = 6\cdot \frac 52 - 6 = 9.
$$
This is confirmed by the data in Table 2.2.

We are also interested in the solutions of the equation
$$
\sum^\infty_{\ell=0} \alpha_\ell(-z)^\ell = -\frac 12.
    \tag6.27
$$
This holds for $z = \lambda_2/5$, because in this case $\lambda_1 = 6$,
and
$$
C_{-\lambda}(x_1) = \sum^\infty_{\ell=0} \alpha_\ell(-\lambda_2/5)^\ell.
$$
But it also holds for $z = \lambda^D_1$, because in this case $\lambda_0
= 6$.  In fact it is easy to see that $\lambda^D_1$ is the smallest
solution of (6.27) (there are infinitely many other choices of $\lambda$
with either $\lambda_1 = 6$ or $\lambda_0 = 6$).  Figure 6.1 shows the
values on $V_1$ of the function $C_{-\lambda}$ in these cases.

\bigskip

\centerline{\BoxedEPSF{: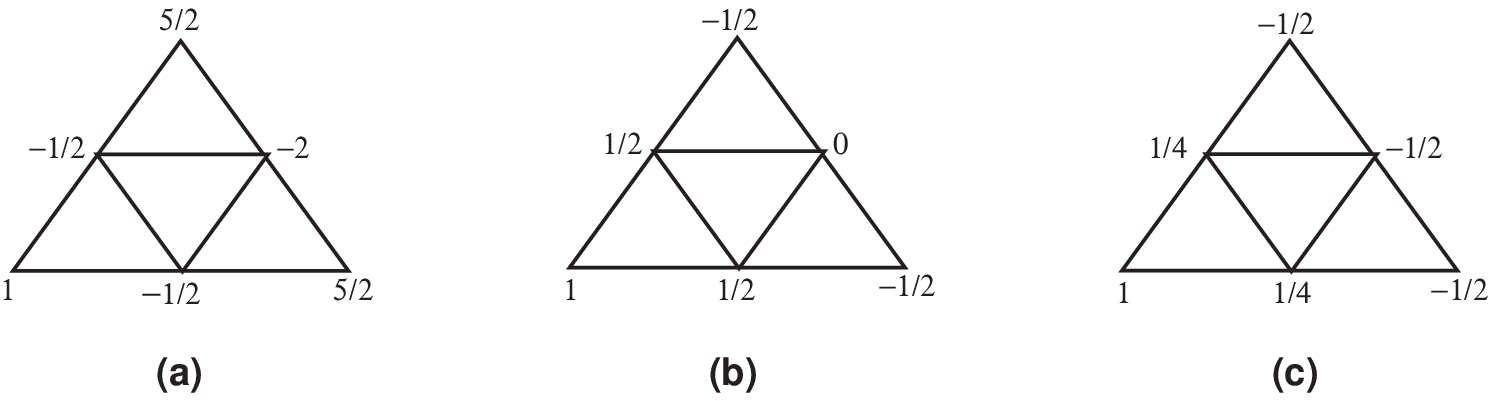 scaled 700}}

\medskip
\noindent
\underbar{Figure 6.1}:  The values of $C_{-\lambda}(x)$ on $V_1$ vertices
for (a) $\lambda_0 = -6$ and $\lambda_1 = 6$, (b)~$\lambda_0 = 6$ and
$\lambda_1 = 2$, (c) $\lambda_0 = 6$ and $\lambda_1 = 3$.


\bigskip
We can now explain why the recursion relation (2.11) for $\beta_j$ is
unstable.  It is clear by inspection that the middle term on the right
side of (2.11) is much larger than the other terms, so we would expect
that a solution of (2.11) would be close to a solution of
$$
\wt\beta_j = -\frac 23 \sum^{j-1}_{\ell=0} \alpha_{j-\ell}
         5^{\ell-j}\wt\beta_\ell,
$$
which may be rewritten as
$$
-\frac 12 = \sum^j_{\ell=0} \alpha_\ell 5^{-\ell}
      \frac{\wt\beta_{j-\ell}}{\wt\beta_j}.
   \tag6.28
$$
If we look for a solution of (6.28) of the form $\wt\beta_j = (-5z)^{-j}$
then we obtain $\ds\sum^j_{\ell=0} \alpha_\ell(-z)^\ell = -\frac 12$,
which is very close to (6.27) in view of the very rapid decay of
$\alpha_\ell$.  The solution to (6.28) should thus be an infinite linear
combination of exponential solutions with $z$ a solution to (6.27).  In
the generic case the dominant term should correspond to the smallest
solution of (6.27).  Thus we expect the solution to (6.28) to behave like
a multiple of $(-5\lambda^D_1)^{-j}$, and numerical computations confirm
this.  This pseudo--solution  of (2.11)  attracts any approximate
solution of (2.11) that strays from the exact solution.

A related observation is that $\ds\sum^\infty_{\ell=0}
\alpha_\ell(-z)^\ell = 1$ holds for $z = \lambda^D_2 \approx
55.885828\ldots$ by (6.15), since in this case $\lambda_0 = 0$ and
$\lambda_1 = 5$.  In the form $\ds\sum^\infty_{\ell=1}
\alpha_\ell(-\lambda^D_2)^\ell = 0$ this suggests that the entries of the
matrix $\sigma(\alpha)^{-1}$, which are just  $6T_j$, should decay like
$(-\lambda^D_2)^{-j}$.  The numerical data in Table 6.3 confirms this.
This explains the instability in the recursion relation for $\{t_j\}$.


$$
\matrix
j   &\quad T_j   &(-\lambda^D_2)^j T_j \\   \vspace{.08in}
0    &\quad 1.    &1.  \\
1    &\quad -.03333333333    &\quad 1.862860915  \\
2    &\quad .0007407407407   &\quad 2.313500526  \\
3    &\quad -.00001433691756    &\quad 2.502423700 \\
4    &\quad .2637965601 \ 10^{-6}    &\quad 2.573213790  \\
5    &\quad -.4766054541 \ 10^{-8}    &\quad 2.598169232  \\
6    &\quad .8556101104 \ 10^{-10}    &\quad 2.606669803  \\
7    &\quad -.1532663873 \ 10^{-11}    &\quad 2.609508520  \\
8    &\quad .2743475872 \ 10^{-13}    &\quad 2.610445492   \\
9    &\quad -.4909650195\  10^{-15}    &\quad 2.610752605   \\
10   &\quad .8785480907 \ 10^{-17}    &\quad 2.610852844   \\
11   &\quad -.1572060595\  10^{-18}   &\quad 2.610885478   \\
12   &\quad .2812997595 \ 10^{-20}    &\quad 2.610896085   \\
13   &\quad -.5033478852 \ 10^{-22}    &\quad 2.610899530  \\
14   &\quad .9006721805\  10^{-24}    &\quad 2.610900647   \\
15   &\quad -.1611629185 \ 10^{-25}    &\quad 2.610901010   \\
16   &\quad .2883788845 \ 10^{-27}    &\quad 2.610901127   \\
17   &\quad -.5160143489  \ 10^{-29}    &\quad 2.610901165   \\
18   &\quad .9233366935 \ 10^{-31}    &\quad 2.610901177   \\
19   &\quad -.1652183992 \ 10^{-32}    &\quad 2.610901182  \\
20   &\quad .2956355963 \ 10^{-34}    &\quad 2.610901182   \\
\endmatrix
$$

\centerline{Table 6.3}
\bigskip
We also observe that the values of $C_{-\lambda_2}(x)$ given in Figure 6.1
(a) show that the rearranged power series at $q_1$ does not converge to
$C_{-\lambda_2}$ outside the cell $F_1(SG)$.  Indeed, the even part of
the power series about $q_1$, if it converged in SG, would have to be
$\frac 52 \sum(-\lambda_2)^j P^{(1)}_{j1}(x)$, which gives the incorrect
value of $25/4$ for $\frac 12(C_{-\lambda_2}(q_0) + C_{-\lambda_2}(q_2))
= 7/4$.

\bye